\documentclass{article}
\usepackage[utf8]{inputenc}
\usepackage[margin =2.4cm]{geometry}
\usepackage{amssymb}
\usepackage{amsmath}
\usepackage{amsthm}
\usepackage{enumerate} 
\usepackage{graphicx}
\usepackage{csquotes}
\usepackage{subcaption}
\usepackage{lmodern}
\usepackage{booktabs}
\usepackage{subcaption}

\newtheorem{theorem}{Theorem}
\newtheorem{proposition}{Proposition}

\newtheorem{lemma}{Lemma}

\title{{\fontsize{25pt}{40pt}\selectfont Evaluating the Sharpness and Limitations of Bounds on the Frobenius Number}}

    \author{} 
    \date{} 

\begin{document}
\maketitle
\vspace{-14mm}
\begin{center}
    {\fontsize{16pt}{40pt}\selectfont Aled Williams} 
    \\  \vspace{4.0mm}
    Department of Mathematics \\
    London School of Economics and Political Science \\
    London, UK \\
    \texttt{a.e.williams1@lse.ac.uk} \\
    \vspace{5mm} 
\end{center}



\begin{abstract}
In this paper we study the (classical) Frobenius problem, namely the problem of finding the largest integer that cannot be represented as a nonnegative integral combination of given relatively prime (strictly) positive integers (known as the Frobenius number). We firstly compare several upper bounds on the Frobenius number, assessing their relative tightness through both theoretical arguments and Monte Carlo simulations. We then explore whether a general upper bound with a worst-case exponent strictly less than quadratic can exist, and formally demonstrate that such an improvement is impossible. These findings offer new insights into the structural properties of established bounds and underscore inherent constraints for future refinement.

\vspace{1.0mm}
\end{abstract}


\section{Introduction}
Let $\boldsymbol{a}$ be a positive integral $n$-dimensional primitive vector, i.e. $\boldsymbol{a} = (a_1, \ldots, a_n)^T \in \mathbb{Z}^n_{>0}$ with $\gcd(\boldsymbol{a}) := \gcd(a_1, \ldots, a_n) =1$. In what follows, we exclude the case $n=1$ and assume that the dimension $n \ge 2$. In particular, without loss of generality, we assume the following conditions:
\begin{equation} \label{conditions on a}
\boldsymbol{a} = (a_1, \ldots, a_n)^T \in \mathbb{Z}^n_{>0}\,, n \ge 2 \text{ and } \gcd(\boldsymbol{a}):= \gcd(a_1, \ldots, a_n) = 1.
\end{equation}
The \textit{Frobenius number} of $\boldsymbol{a}$, denoted by $F(\boldsymbol{a})$, is the largest integer that cannot be represented as a nonnegative integral combination of the $a_i$'s, i.e.
$$
F(\boldsymbol{a}) := \max \left\{ b \in \mathbb{Z} : b \ne \boldsymbol{a}^T \boldsymbol{z} \text{ for all } \boldsymbol{z} \in \mathbb{Z}^n_{\ge 0} \right\} ,
$$
where $\boldsymbol{a}^T$ denotes the transpose of $\boldsymbol{a}$. Note for completeness that the Frobenius problem, namely the problem of finding the Frobenius number $F(\boldsymbol{a})$ given $\boldsymbol{a}$, is also known by other names within the literature including the money-changing problem (or the money-changing problem of Frobenius, or the coin-exchange problem of Frobenius) \cite{wilf1978circle, tripathi2003variation, bocker2005money}, the coin problem (or the Frobenius coin problem) \cite{boju2007math, spivey2007quadratic} and the Diophantine problem of Frobenius \cite{selmer1977linear, rodseth1978linear}. From a geometric viewpoint, $F(\boldsymbol{a})$ is the maximal right-hand side such that the \textit{knapsack polytope} 
$$
P(\boldsymbol{a},b) = \left\{ \boldsymbol{x} \in \mathbb{R}^n_{\ge 0} : \boldsymbol{a}^T \boldsymbol{x} = b
\right\}
$$
does not contain integral points. It should be noted that \eqref{conditions on a} are indeed necessary and sufficient conditions for the existence of the Frobenius number. 

\vspace{2.0mm}


Note that instead of the conditions \eqref{conditions on a}, some authors instead assume the stronger condition that all the entries of the vector are pairwise coprime, i.e. 
\begin{equation} \label{stronger conditions on a}
\boldsymbol{a} = (a_1, \ldots, a_n)^T \in \mathbb{Z}^n_{>0}\,, n \ge 2 \text{ and } \gcd(a_i, a_j) = 1 \text{ for any } i,j \in \{1,2,\ldots, n\} \text{ with } i \ne j.
\end{equation}
It should be noted that not all integral vectors satisfying \eqref{conditions on a} also satisfy the stronger conditions \eqref{stronger conditions on a}. For example, the vector $\boldsymbol{a} = (6,10,15)^T$ satisfies \eqref{conditions on a} but does not satisfy \eqref{stronger conditions on a}.

\vspace{2.0mm}

There is a very rich history on Frobenius numbers and the book \cite{alfonsin2005diophantine} provides a very good survey of the problem. Note that computing the Frobenius number in general is $\mathcal{NP}$-hard \cite{ramirez1996complexity} (which was proved via a reduction to the integer knapsack problem), however, if the number of integers $n$ is fixed, then a polynomial time algorithm to calculate $F(\boldsymbol{a})$ exists \cite{kannan1992lattice}. 

\vspace{2.0mm}

If $n=2$, it is well-known (most likely due to Sylvester \cite{sylvester1884problem}) that 
\begin{equation} \label{Sylvester 2nd Frobenius bound}
\begin{aligned}
F(a_1, a_2) &= a_1 a_2 - a_1 - a_2 \\
&= (a_1 - 1) (a_2 - 1) - 1\,.
\end{aligned}
\end{equation}
In contrast to the case when $n=2$, it was shown by Curtis \cite{curtis1990formulas} that no general closed formula exists for the Frobenius number if $n>2$. In light of this, there has been a great deal of research into producing upper bounds on $F(\boldsymbol{a})$. These bounds share the property that in the worst-case they are of quadratic order with respect to the maximum absolute valued entry of $\boldsymbol{a}$, which will be denoted by $\| \boldsymbol{a} \|_{\infty} = \max_i |a_i|$. Further, let $\| \cdot \|_{2}$ denote the $\ell_2$-norm (or Euclidean norm).
In particular, upon assuming that \eqref{conditions on a} and $a_1 \le a_2 \le \cdots \le a_n$ hold, such bounds include the classical bound by Erd{\H o}s and Graham \cite[Theorem 1]{erdos1972linear}
\begin{equation*}
F(\boldsymbol{a}) \le 2 a_{n-1} \left\lfloor \frac{a_n}{n} \right\rfloor - a_n ,
\end{equation*}
by Schur (according to Brauer \cite{brauer1942problem})
\begin{equation*}
F(\boldsymbol{a}) \le (a_1 - 1) (a_n - 1) - 1,
\end{equation*} 
by Selmer \cite{selmer1977linear} (for $a_1 \ge n$)
\begin{equation*} 
F(\boldsymbol{a}) \le 2 a_n \left\lfloor \frac{a_1}{n} \right\rfloor - a_1 ,
\end{equation*}
by Vitek \cite[Theorem 5]{vitek1975bounds} (for $n \ge 3$)
\begin{equation*}
F(\boldsymbol{a}) \le \frac{1}{2}(a_2 - 1) (a_n - 2) - 1,
\end{equation*}
and by Fukshansky and Robins \cite[Equation 29]{fukshansky2007frobenius}
\begin{equation*}
F(\boldsymbol{a}) \le \left\lfloor \frac{(n-1)^{2} \,\Gamma (\frac{n+1}{2})}
{\pi^{(n-1)/ 2}} \sum_{i=1}^{n} a_{i} \sqrt{\|\boldsymbol{a}\|_{2}^{2}-a_{i}^{2}}+1 \right\rfloor ,
\end{equation*}
where $\Gamma (\cdot)$ and $\lfloor \, \cdot \, \rfloor$ denote Euler's gamma and the standard floor functions, respectively. 

\vspace{2.0mm}

It should be noted that the bound of Selmer \cite{selmer1977linear} 
\begin{equation} \label{Selmer_bound_n}
F(\boldsymbol{a}) \le 2 a_n \left\lfloor \frac{a_1}{n} \right\rfloor - a_1,
\end{equation}
with $a_1 \ge n$, is widely referenced in books and papers, but is frequently misstated. In particular, many sources give insufficient attention to the underlying assumptions on the vector $\boldsymbol{a}$. Notably, this upper bound does not necessarily hold under the weaker conditions \eqref{conditions on a}. 

\vspace{2.0mm}

\begin{proposition} \label{prop1}
The upper bound \eqref{Selmer_bound_n} of Selmer \cite{selmer1977linear} does not necessarily hold unless the stronger conditions \eqref{stronger conditions on a} hold. This requirement remains even if the weaker conditions \eqref{conditions on a} are met. 
\end{proposition}

\vspace{2.0mm}

Recall that the above bounds share the property that in the worst-case they are of quadratic order and depend on the maximum entry of $\boldsymbol{a}$. It was shown by Williams and Haijima \cite{williams2023considering} that any upper bound on the Frobenius number $F(\boldsymbol{a})$ must inherently depend on $\| \boldsymbol{a} \|_{\infty}$ if the integral vector $\boldsymbol{a}$ satisfies only the (weaker) conditions \eqref{conditions on a}. Later in this paper (in Section 4), we investigate whether a general upper bound on $F(\boldsymbol{a})$ can be obtained under these weaker conditions with a worst-case exponent strictly lower than quadratic.

\vspace{2.0mm}

If the entries of $\boldsymbol{a}$ satisfy the stronger conditions \eqref{stronger conditions on a} that they are pairwise coprime, then, following an argument closely related to Beck et al. \cite{beck2002frobenius}, Williams and Haijima \cite{williams2023considering} establish the bound 
\begin{equation} \label{Williams_corrected}
F(\boldsymbol{a}) \le \frac{1}{2}\left(\sqrt{\frac{1}{3}\left(a_1+a_2+a_3\right)\left(a_1+a_2+a_3+2 a_1 a_2 a_3\right)+\frac{8}{3}\left(a_1 a_2+a_2 a_3+a_3 a_1\right)}-a_1-a_2-a_3\right) . 
\end{equation}
Notably, the original proof of Beck et al. \cite{beck2002frobenius} contained a subtle error, which altered the value of their initially claimed upper bound of
\begin{equation} \label{Beck_et_al_bound}
F(\boldsymbol{a}) \le \frac{1}{2} \left( \sqrt{a_1 a_2 a_3 \left( a_1 + a_2 + a_3 \right)} - a_1 - a_2 - a_3 \right).
\end{equation}
It was shown by Williams and Haijima 
\cite{williams2023considering}
that despite the subtle error, the original bound \eqref{Beck_et_al_bound} remains valid assuming the stronger conditions \eqref{stronger conditions on a} hold, although it turns out that \eqref{Williams_corrected} is tighter than \eqref{Beck_et_al_bound} in all but a relatively \enquote{small} (finite) number of cases. 
Furthermore, it was shown by Williams and Haijima \cite{williams2023considering} that provided the stronger conditions \eqref{stronger conditions on a} hold, then for any $i,j \in \{1,2,\ldots,n\}$ with $i \ne j$ we have 
$$
F(\boldsymbol{a}) \le (a_i - 1)(a_j - 1) - 1.
$$
It should be noted that the above bound tells us that the well-known result \eqref{Sylvester 2nd Frobenius bound} of Sylvester \cite{sylvester1884problem} naturally extends to provide an upper bound for the Frobenius number $F(\boldsymbol{a})$ under the (stronger) conditions \eqref{stronger conditions on a}, i.e. that the entries of $\boldsymbol{a}$ are pairwise coprime. 

\vspace{2.0mm}


\section{Comparison of Bounds Under the GCD Conditions}
In this section, we compare bounds that are valid under the (weaker) conditions \eqref{conditions on a} on $\boldsymbol{a}$. In particular, we examine the bounds by Erd{\H o}s and Graham \cite[Theorem 1]{erdos1972linear}
\begin{equation} \label{erdos_upper}
F(\boldsymbol{a}) \le 2 a_{n-1} \left\lfloor \frac{a_n}{n} \right\rfloor - a_n ,
\end{equation}
by Schur (according to Brauer \cite{brauer1942problem})
\begin{equation} \label{schur_upper}
F(\boldsymbol{a}) \le (a_1 - 1) (a_n - 1) - 1,
\end{equation}
by Vitek \cite[Theorem 5]{vitek1975bounds} (for $n \ge 3$)
\begin{equation} \label{vitek_upper}
F(\boldsymbol{a}) \le \frac{1}{2}(a_2 - 1) (a_n - 2) - 1,
\end{equation}
and by Fukshansky and Robins \cite[Equation 29]{fukshansky2007frobenius}
\begin{equation} \label{Fukshansky_upper}
F(\boldsymbol{a}) \le \left\lfloor \frac{(n-1)^{2} \,\Gamma (\frac{n+1}{2})}
{\pi^{(n-1)/ 2}} \sum_{i=1}^{n} a_{i} \sqrt{\|\boldsymbol{a}\|_{2}^{2}-a_{i}^{2}}+1 \right\rfloor ,
\end{equation}
where $\Gamma (\cdot)$ and $\lfloor \, \cdot \, \rfloor$ denote Euler's gamma and the standard floor functions, respectively. 

\vspace{2.0mm}

It should be noted that in several of the forthcoming figures, the names of above bounds are abbreviated using internal variable names. In particular, \texttt{diffErdos} refers to the bound of Erd{\H o}s and Graham \eqref{erdos_upper}, \texttt{diffSchur} refers to the bound of Schur \eqref{schur_upper}, \texttt{diffVitek} refers to the bound of \eqref{vitek_upper}, and \texttt{diffFukRob} refers to the bound of Fukshansky and Robins \eqref{Fukshansky_upper}.

\vspace{2.0mm}


We begin our analysis by comparing the distribution of errors associated with each bound across various dimensions $n$, under differing constraints on the maximum absolute entry of $\boldsymbol{a}$. In this simulation, we generate integer vectors $\boldsymbol{a}$ satisfying the weaker condition \eqref{conditions on a}, with entries ordered such that $a_1 \le a_2 \le \cdots \le a_n$. For each vector, we compute the Frobenius number $F(\boldsymbol{a})$ using the algorithm of Wilf \cite{wilf1978circle} or Nijenhuis \cite{nijenhuis1979minimal}, along with the corresponding values of the upper bounds under study. The integer vector entries are sampled uniformly at random from the interval $[k, m]$, where $k \ge n$ and $m \in \{ 100, 1000, \, 10,000 \}$, meaning that $\|\boldsymbol{a}\|_\infty \le m$ holds. This process is repeated 100,000 times for each dimension $n$ under consideration.

\vspace{2.0mm}

\begin{figure}[ht!]
    \centering
    \includegraphics[width=1\linewidth]{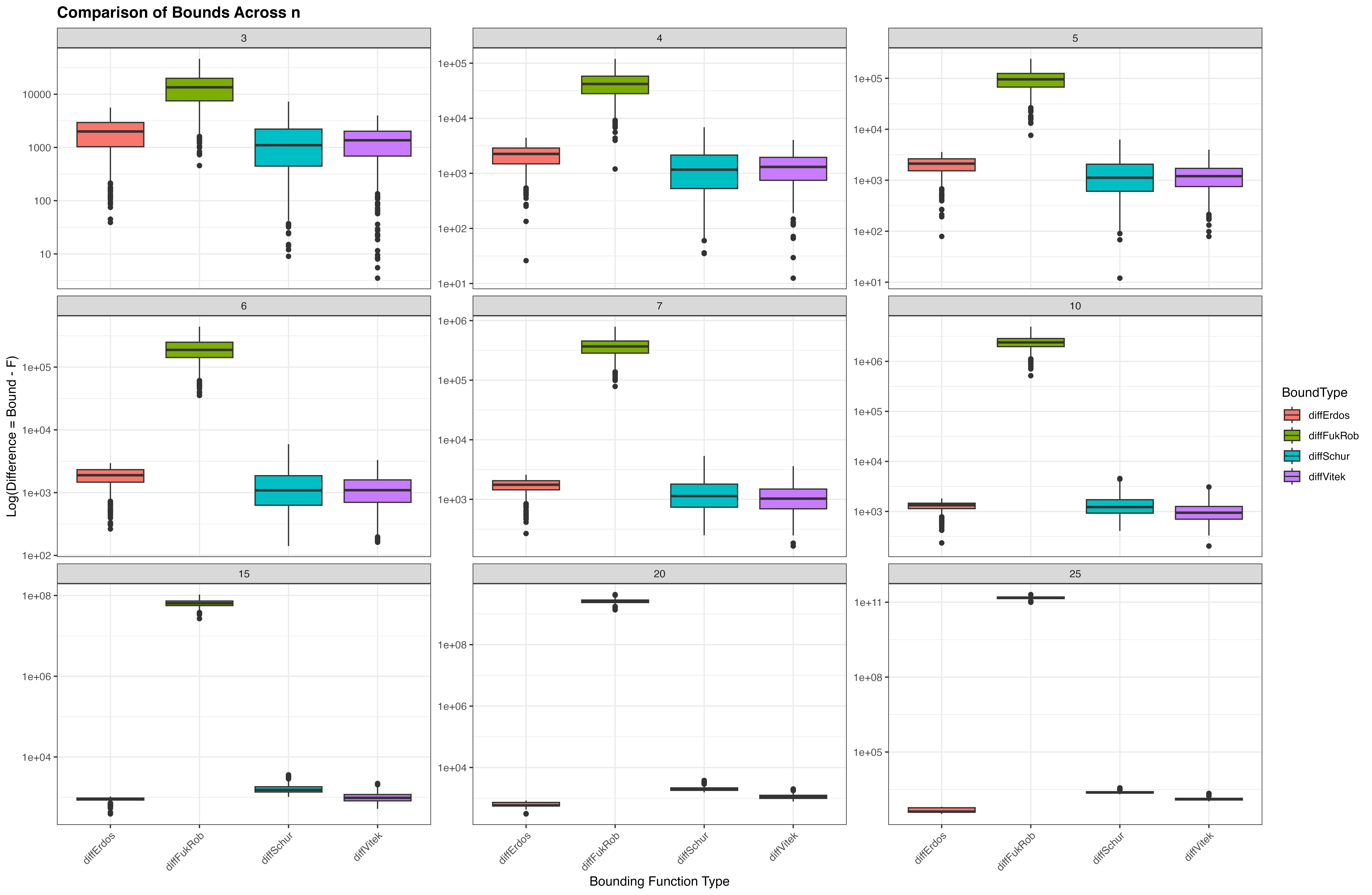}
    \caption{Box plots of the difference, i.e. \(\text{Bound} - F(\boldsymbol{a})\), across dimensions \(n\) with \(\|\boldsymbol{a}\|_\infty \le 100\). The log-scaled vertical axis reveals the comparative tightness and spread of different bounds at \enquote{moderate} input size.}
    \label{fig:boxplot_main_m100}
\end{figure}

Figure \ref{fig:boxplot_main_m100} displays log-scaled box plots of bound errors under a \enquote{moderate} upper bound on the vector entries, namely with $m = 100$. Each box summarise the error distribution, where the centre line indicates the median, and the boxes themselves represent the interquartile range. 
For larger values of $m$, the growth of the error, especially for the Fukshansky and Robins bound, becomes more pronounced. The corresponding box plots for $m = 1000$ and $m = 10000$ are presented in Appendix~\ref{appendix:larger_m}. Observe that for small dimension $n$, the difference is relatively modest for all bounds bar the bound of Fukshansky and Robins \eqref{Fukshansky_upper}. The bounds of Erd{\H o}s and Graham, Schur and Vitek show more stability upon increases in $n$. 

\vspace{2.0mm}

\vspace{2.0mm}

\begin{figure}[ht!]
    \centering
    \includegraphics[width=1\linewidth]{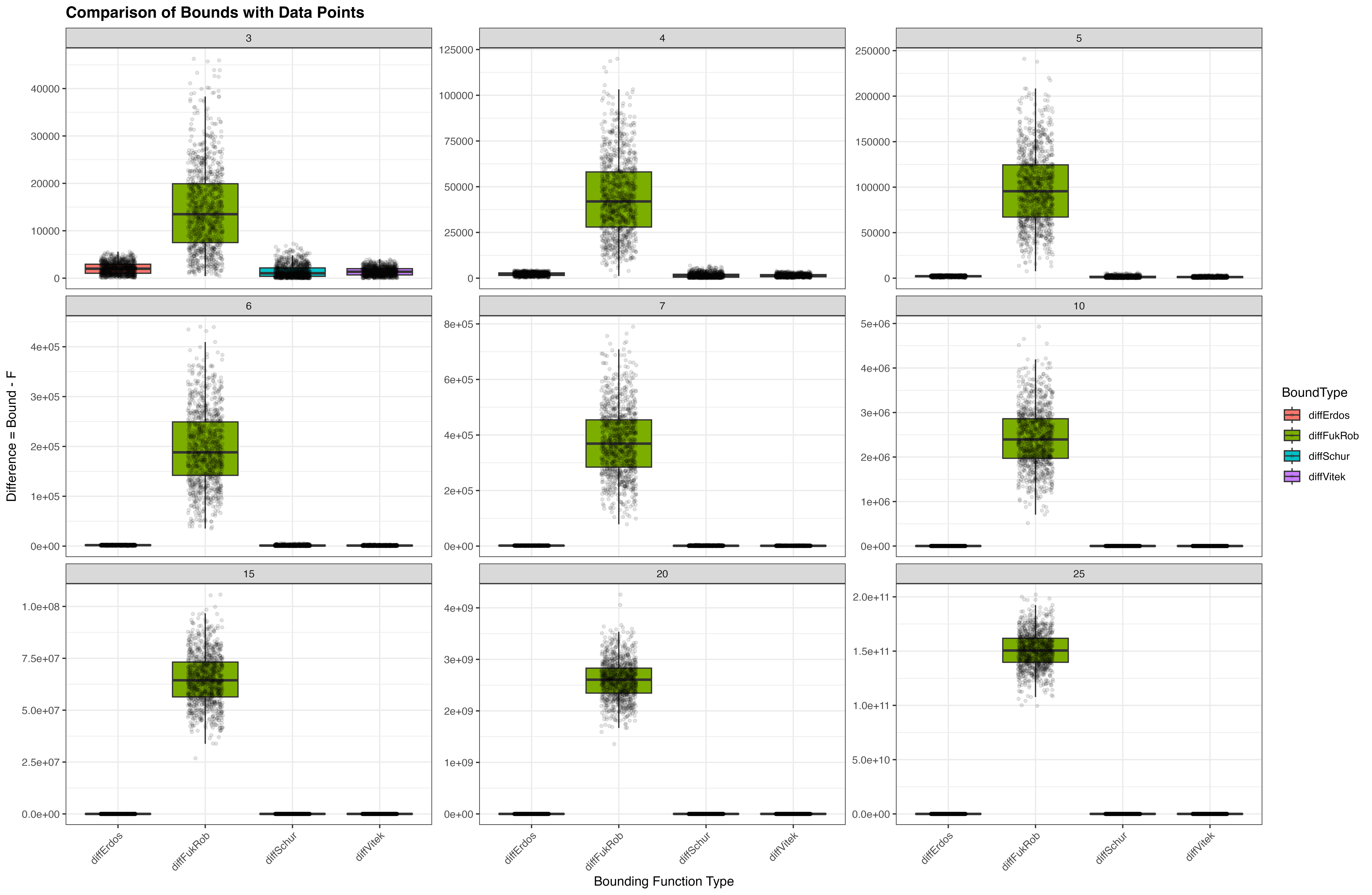}
    \caption{Box plots of the difference, i.e. \(\text{Bound} - F(\boldsymbol{a})\), with jittered sample points with \(\|\boldsymbol{a}\|_\infty \le 100\). Overlaid data points reveal the distributional density and outliers. The log-scaled vertical axis highlights the growing disparity across bounds, especially for the Fukshansky and Robins bound.}
    \label{fig:boxplot_data_m100}
\end{figure}

While the box plots reveal broad differences in the behaviour of the bounds, they obscure the underlying sample density. To address this, Figure \ref{fig:boxplot_data_m100} overlays the raw Monte Carlo data points, making it easier to identify the consistency (or lack thereof) of each bound. Notably, the bounds of Schur, Vitek, and Erdős and Graham display much tighter concentration near the median, while the bound of Fukshansky and Robins exhibits wide dispersion and many outliers. Plots for larger values of $m$ (i.e. 1000 and 10000) are presented in Appendix~\ref{appendix:densities}.

\vspace{2.0mm}

\vspace{2.0mm}

\begin{figure}[ht!]
    \centering
    \includegraphics[width=1\linewidth]{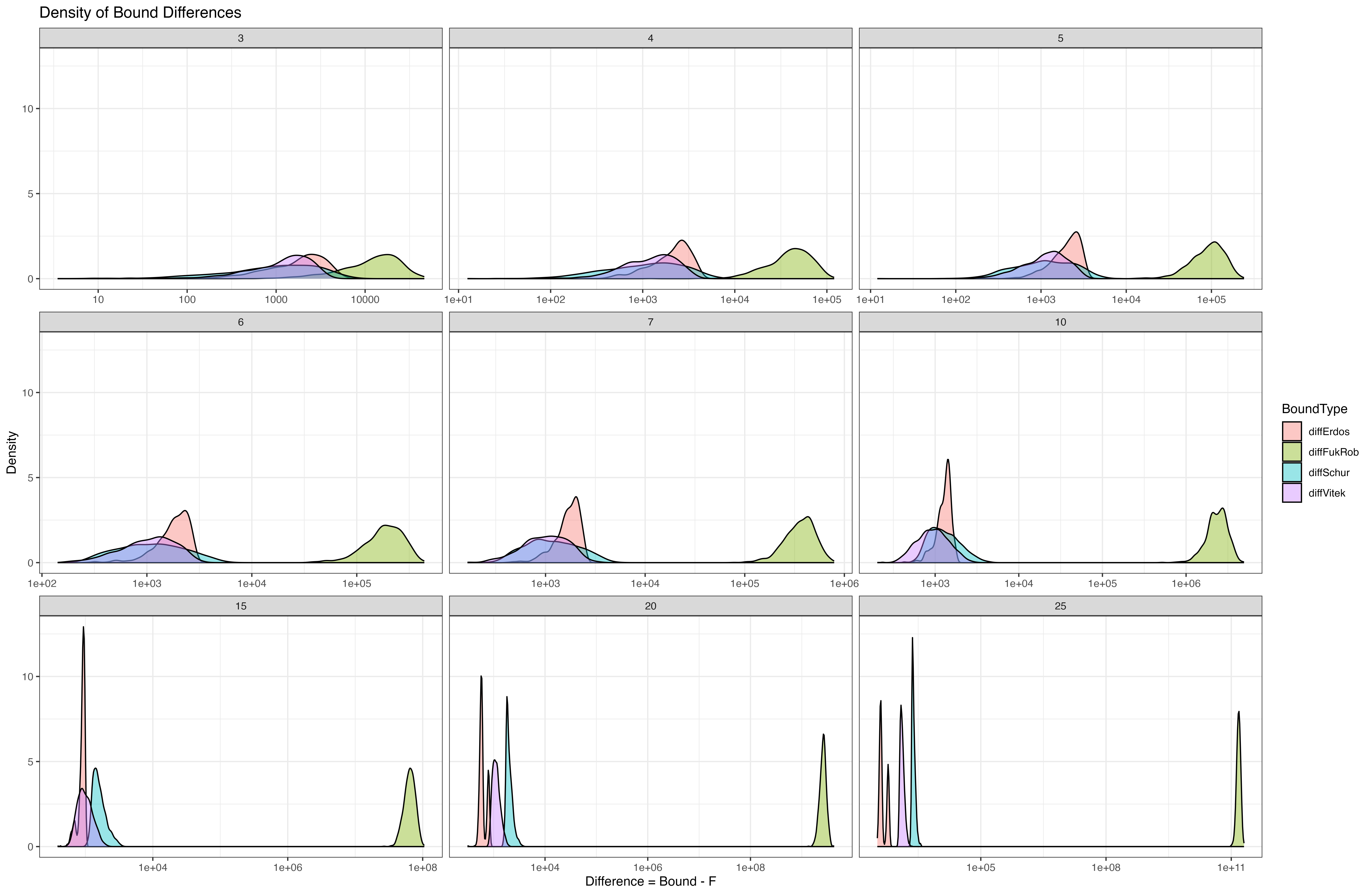}
    \caption{Density plots of the difference, i.e. \(\text{Bound} - F(\boldsymbol{a})\), for \(\|\boldsymbol{a}\|_\infty \le 100\), across different dimensions \(n\). Each bound exhibits a concentrated distribution around its median, with the Fukshansky and Robins bound showing heavier tails indicative of large outliers.}
    \label{fig:density_m100}
\end{figure}

To explore this variation in more detail, Figure~\ref{fig:density_m100} presents density plots of the error distributions for each upper bound. These reveal that while most bounds have unimodal, compact distributions, the Fukshansky and Robins bound exhibits longer tails, confirming its extreme overestimation in some cases, especially as the dimension grows. Density plots for larger values of $m$ are presented in Appendix~\ref{appendix:density_plots}, where the heavy-tailed behaviour of the Fukshansky and Robins bound becomes more pronounced.

\vspace{2.0mm}

\vspace{2.0mm}

\begin{figure}[ht!]
    \centering
    \includegraphics[width=1\linewidth]{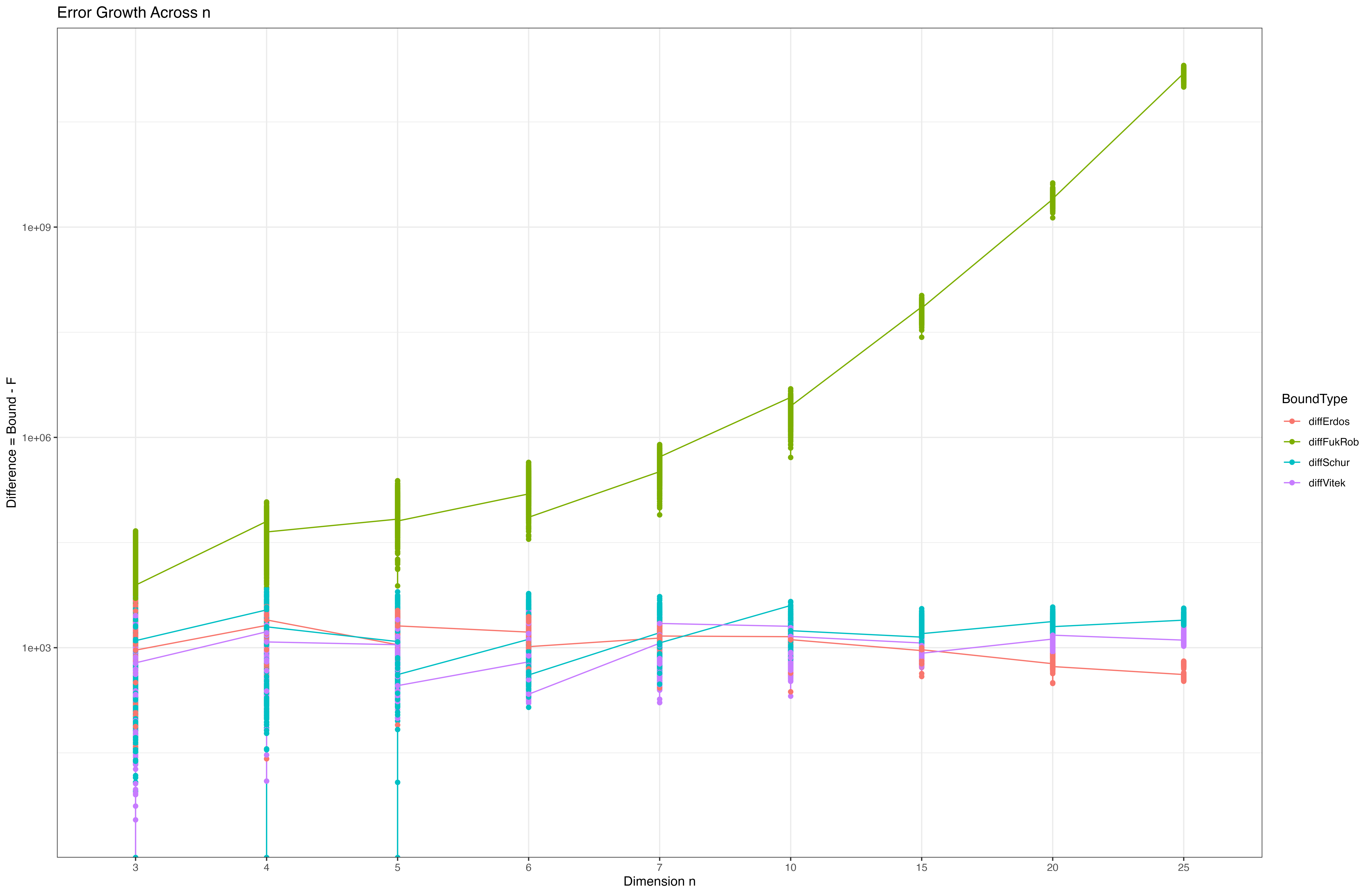}
    \caption{The average difference, i.e. \(\text{Bound} - F(\boldsymbol{a})\), plotted against vector dimension \(n\), with \(\|\boldsymbol{a}\|_\infty \le 100\). The classical bounds exhibit sub-exponential growth, while the Fukshansky and Robins bound grows significantly faster with increasing dimension.}
    \label{fig:line_error_m100}
\end{figure}

Next, we investigate how error magnitudes scale with dimension $n$. Figure~\ref{fig:line_error_m100} plots the average difference, i.e. $\text{Bound} - F(\boldsymbol{a})$, for each bounding function against the dimension $n$. While all bounds show increasing error, their growth rates differ significantly. Classical bounds such as Schur, Vitek, and Erdős grow modestly, consistent with their polynomial dependence on $\|\boldsymbol{a}\|_\infty$, whereas the Fukshansky and Robins bound increases at a much quicker rate. To maintain clarity of presentation, results for $\|\boldsymbol{a}\|_\infty \le 100$ are presented in the main text, while larger-scale comparisons are included in Appendix~\ref{appendix:line_error_plots}.

\vspace{2.0mm}

\vspace{2.0mm}

\begin{figure}[ht!]
    \centering
    \includegraphics[width=1\linewidth]{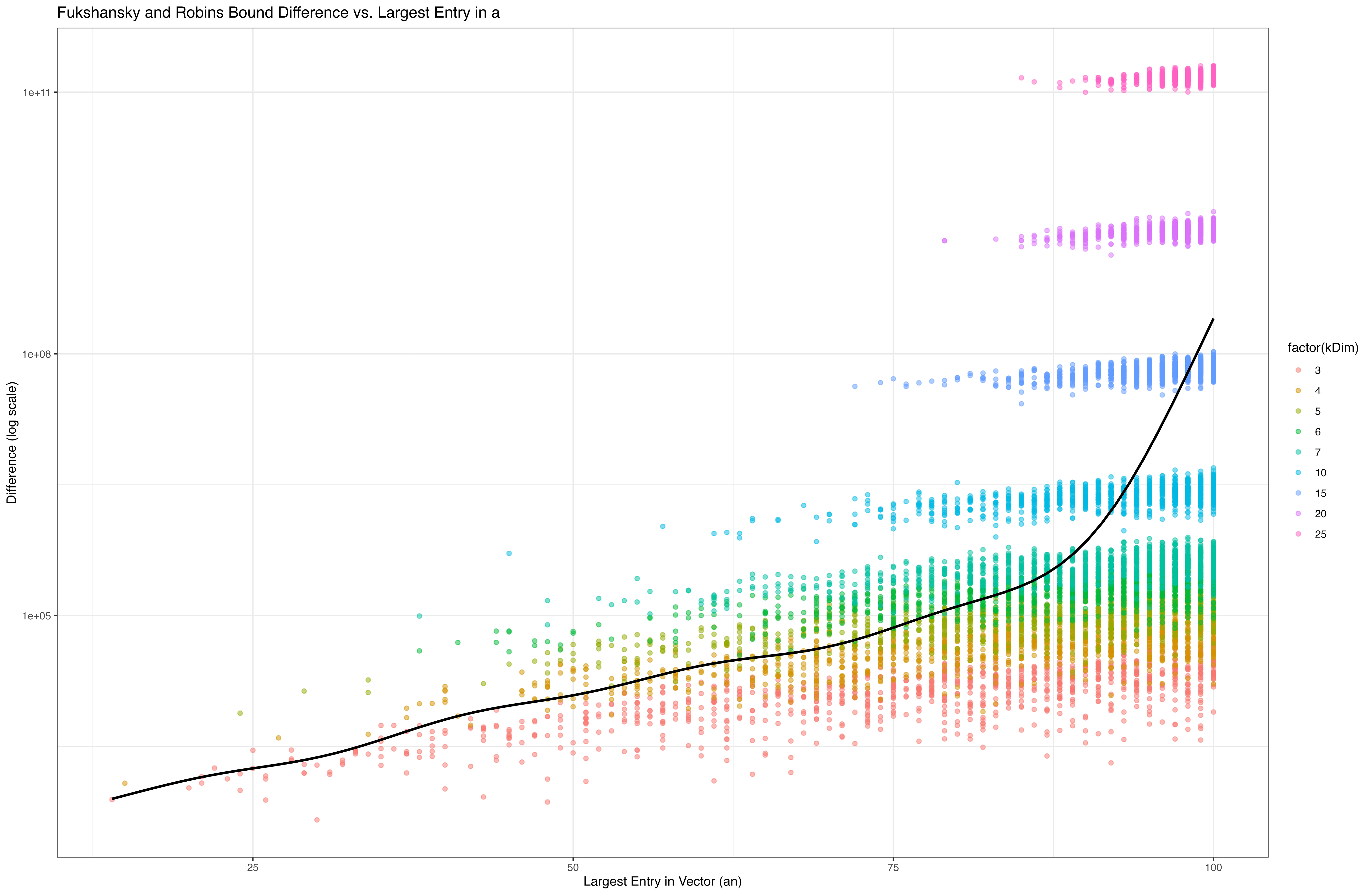}
    \caption{Scatter plot of the difference \(\text{FukRob} - F(\boldsymbol{a})\) plotted against the maximum absolute entry \(\|\boldsymbol{a}\|_{\infty}\) with \(\|\boldsymbol{a}\|_\infty \le 100\). The colours denote the dimension \(n\). The log-scaled vertical axis highlights the exponential relationship between bound error and the size of the input coefficients.}
    \label{fig:scatter_fukrob_m100}
\end{figure}

Finally, we explore how the error in the Fukshansky and Robins bound relates to the maximum entry in the vector $\boldsymbol{a}$. Figure~\ref{fig:scatter_fukrob_m100} displays a scatter plot of the difference $\text{FukRob} - F(\boldsymbol{a})$ against the maximum absolute valued entry $\| \boldsymbol{a} \|_{\infty}$, where points are coloured by their dimension $n$. A non-linear trend emerges, namely as $\|\boldsymbol{a}\|_\infty$ increases, the bound error grows rapidly. This is even more pronounced for higher dimensions.
Results for larger bounds on $\|\boldsymbol{a}\|_\infty$ are presented in Appendix~\ref{appendix:scatter_plots}.

\vspace{2.0mm}

\vspace{2.0mm}

The empirical results in Figures \ref{fig:boxplot_main_m100}-\ref{fig:scatter_fukrob_m100} consistently show that the Fukshansky and Robins bound grows much more rapidly than the classical bounds of Schur, Erdős and Graham, or Vitek, particularly as either the dimension $n$ or the maximum entry $\|\boldsymbol{a}\|_\infty$ increases. 
Despite the fact that all these bounds are valid upper bounds assuming the (weaker) conditions \eqref{conditions on a}, their \enquote{internal structure} determines how tightly they bound $F(\boldsymbol{a})$ in practice. 
In particular, the Fukshansky and Robins bound includes both an $\ell_2$-norm term and a constant involving the gamma function, which together introduce super-polynomial growth which is not present in the classical and more algebraically simple bounds. To better understand why the Fukshansky and Robins bound exhibits such scaling behaviour, we now carefully consider its algebraic structure and compare it with the asymptotic forms of the classical bounds.

\vspace{2.0mm}

Recall that the Fukshansky and Robins upper bound is given by 
$$
\left\lfloor \frac{(n-1)^{2} \,\Gamma (\frac{n+1}{2})}
{\pi^{(n-1)/ 2}} \sum_{i=1}^{n} a_{i} \sqrt{\|\boldsymbol{a}\|_{2}^{2}-a_{i}^{2}}+1 \right\rfloor.
$$
We firstly focus on the factor
$$
C(n) = \frac{(n-1)^2\, \Gamma\!\left(\frac{n+1}{2}\right)}{\pi^{\frac{n-1}{2}}},
$$
which depends on the dimension $n$. Using Stirling’s approximation for the gamma function (see e.g. \cite[Chapter 2]{feller1991introduction}), for large $n$, we have 
$$
\Gamma\!\left(\frac{n+1}{2}\right) \sim \sqrt{2\pi} \left(\frac{n-1}{2}\right)^{\frac{n-1}{2}} e^{-\frac{n-1}{2}}, 
$$
where $\sim$ denotes asymptotic equivalence (i.e. that the ratio of the two expressions tends to 1 as some parameter tends to infinity). Thus, we have 
$$
C(n) \sim (n-1)^2 \sqrt{2\pi} \left(\frac{n-1}{2\pi}\right)^{\frac{n-1}{2}} e^{-\frac{n-1}{2}}
$$
for large $n$.

\vspace{2.0mm}

\begin{figure}[ht!]
    \centering
    \includegraphics[width=0.8\linewidth]{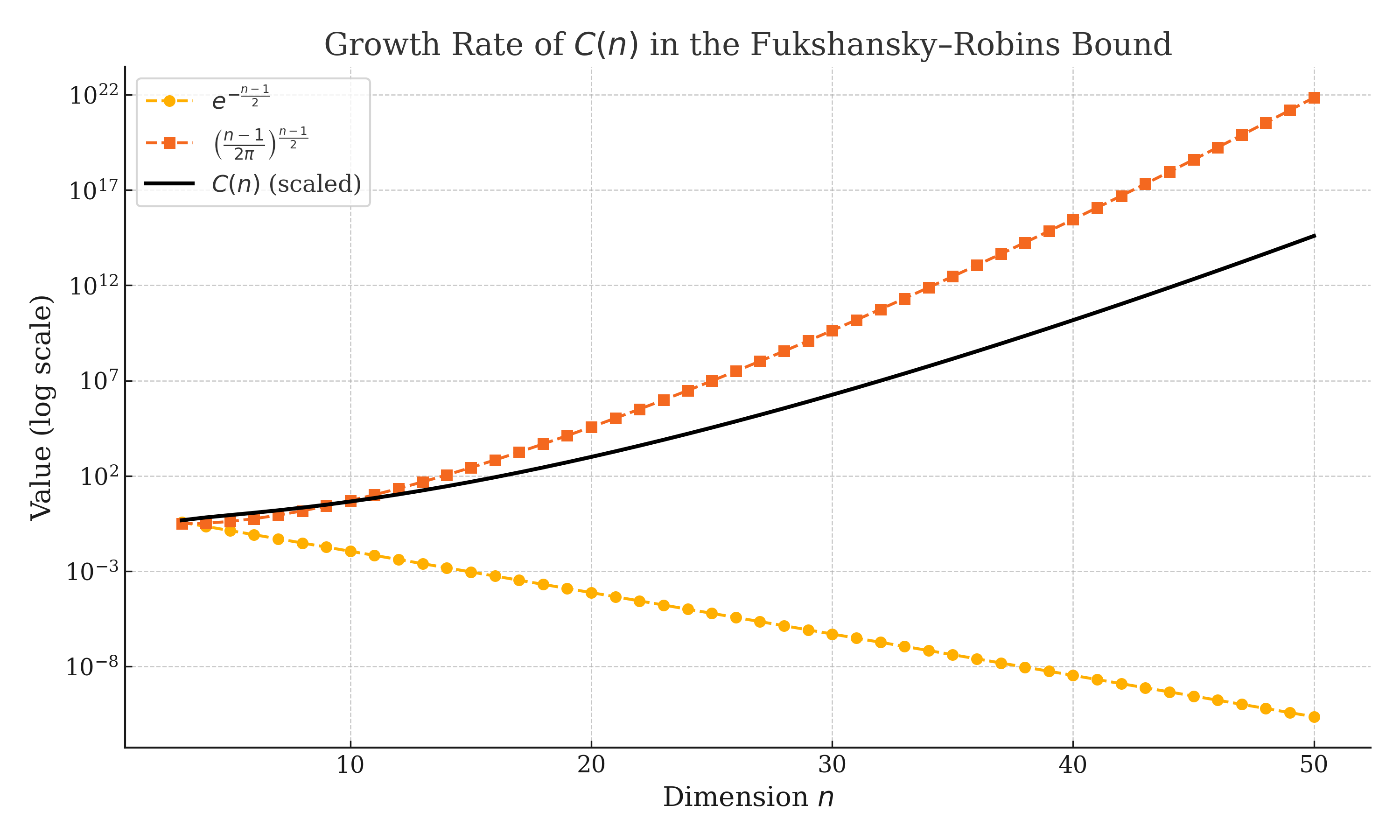}
    \caption{Growth behaviour of the leading constant \(C(n)\) in the Fukshansky and Robins bound. Although the term \(e^{-(n-1)/2}\) decays exponentially, the dominating term \(\left(\frac{n-1}{2\pi}\right)^{(n-1)/2}\) grows super-exponentially, leading to rapid overall growth of \(C(n)\). This explains the empirical divergence of the bound in higher dimensions.}
    \label{fig:fukrob_Cn_growth}
\end{figure}

Observe that the exponential term $e^{-\frac{n-1}{2}}$ decays in $n$, while the term $(\frac{n-1}{2\pi})^{\frac{n-1}{2}}$ grows super-exponentially for large $n$. Thus, it follows that $C(n)$ grows quickly with increases in $n$. To better visualise this, Figure~\ref{fig:fukrob_Cn_growth} plots on a log-scale the growth of the component terms of $C(n)$ and a scaled $C(n)$ separately. While the exponential decay $e^{-(n-1)/2}$ rapidly decreases, it is ultimately overwhelmed by the super-exponential growth of $\left( \frac{n-1}{2\pi} \right)^{(n-1)/2}$. Their product $C(n)$ increases rapidly with the dimension $n$.

\vspace{2.0mm}

\vspace{2.0mm}

Next, consider the sum 
$$
S(\boldsymbol{a}) = \sum_{i=1}^{n} a_i\, \sqrt{\|\boldsymbol{a}\|_2^2 - a_i^2}. 
$$
If for simplicity we assume that all $a_i$'s are of the same order (i.e. around some typical value $A = \Theta(m)$ for some $m \in \mathbb{Z}_{>0}$), then
$$
\|\boldsymbol{a}\|_2^2 \approx n \, A^2, 
$$
and each term of the sum is approximated by 
$$
a_i\, \sqrt{n A^2 - A^2} \approx A\, \sqrt{(n-1)A^2} = A^2 \sqrt{n-1}.
$$
Since there are $n$ such terms, we yield
$$
S(\boldsymbol{a}) \approx n A^2 \sqrt{n-1}
$$
when all $a_i$'s are of the same order. Thus, the Fukshansky and Robins bound is roughly of order $C(n) \cdot S(\boldsymbol{a})$ and grows extremely quickly with $n$, even for fixed $A$ (or $m$). Most classical bounds, in contrast, namely those of Schur, Erdős and Graham, and Vitek, are of at most quadratic order in the $a_i$'s, and therefore scale approximately like $\Theta(A^2)$ when the components of the vector $\boldsymbol{a}$ are uniformly bounded. This fundamental difference in asymptotic behaviour explains the rapidly growing overestimation observed empirically in the Fukshansky and Robins bound.

\vspace{2.0mm}

In summary, these visualisations paint a consistent picture, namely that while \eqref{erdos_upper}, \eqref{schur_upper}, \eqref{vitek_upper} and \eqref{Fukshansky_upper} are valid upper bounds, their practical sharpness varies significantly. The Schur, Erdős and Graham, and Vitek bounds maintain relatively stable performance across dimensions and input sizes, making them somewhat more suitable for practical estimation. The Fukshansky and Robins bound, in contrast, though geometrically motivated, exhibits substantial overestimation, especially as either the dimension or the maximum entry size increases. This behaviour is consistent with its theoretical structure, which combines gamma functions and $\ell_2$-norm expressions, both of which scale unfavourably in high-dimensional settings.

\vspace{2.0mm}


\section{Comparison of Bounds Under the Pairwise Coprime-Conditions}
In this section, we compare the bounds by Selmer \cite{selmer1977linear} (for $a_1 \ge n$)
\begin{equation} \label{Selmer_bound_n_1}
F(\boldsymbol{a}) \le 2 a_n \left\lfloor \frac{a_1}{n} \right\rfloor - a_1,
\end{equation}
by Beck et al. \cite{beck2002frobenius}
\begin{equation} \label{Beck_et_al_bound_1}
F(\boldsymbol{a}) \le \frac{1}{2} \left( \sqrt{a_1 a_2 a_3 \left( a_1 + a_2 + a_3 \right)} - a_1 - a_2 - a_3 \right), 
\end{equation}
the first bound by Williams and Haijima \cite{williams2023considering}
\begin{equation} \label{Williams_corrected_1}
F(\boldsymbol{a}) \le \frac{1}{2}\left(\sqrt{\frac{1}{3}\left(a_1+a_2+a_3\right)\left(a_1+a_2+a_3+2 a_1 a_2 a_3\right)+\frac{8}{3}\left(a_1 a_2+a_2 a_3+a_3 a_1\right)}-a_1-a_2-a_3\right),
\end{equation}
and the second bound by Williams and Haijima \cite{williams2023considering}
\begin{equation} \label{williams_sylvester_extension}
F(\boldsymbol{a}) \le \min \left\{ (a_i - 1)(a_j - 1) - 1 : i, j \in \{1,2,\ldots, n\} \text{ with } i \ne j \right\},
\end{equation}
which generalises the well-known result \eqref{Sylvester 2nd Frobenius bound} of Sylvester \cite{sylvester1884problem}.

\vspace{2.0mm}

It should be noted that in several of the forthcoming figures, the names of bounds are similarly abbreviated using internal variable names. In particular, \texttt{diffSelmer} refers to the bound of \eqref{Selmer_bound_n_1}, \texttt{diffBeck} refers to the bound of Beck et al. \eqref{Beck_et_al_bound_1}, \texttt{diffWHCorr} refers to the first bound of Williams and Haijima \eqref{Williams_corrected_1}, and \texttt{diffWHMinSyl} refers to the second bound of Williams and Haijima \eqref{williams_sylvester_extension}. 

\vspace{2.0mm}

We begin by empirically evaluating the difference between each bound and the true Frobenius number $F(\boldsymbol{a})$ for randomly generated integer vectors $\boldsymbol{a}$ across various dimensions $n$. In this simulation, we generate integer vectors $\boldsymbol{a}$ satisfying the (stronger) pairwise coprimality conditions \eqref{stronger conditions on a}, with entries ordered such that $a_1 \le a_2 \le \cdots \le a_n$. The integer entries $a_i$ are sampled uniformly at random from the interval $[k, m]$, where $k \ge n$ and $m \in \{ 10,000, \, 100,000 \}$, meaning that $\|\boldsymbol{a}\|_\infty \le m$ holds. This process is repeated 100,000 times for each dimension $n$ under consideration. Note that we set $m=100,000$ in the following, where corresponding plots for $m=10,000$ are presented in appendices. 

\vspace{2.0mm}

\begin{figure}[ht!]
\centering
\includegraphics[width=1\linewidth]{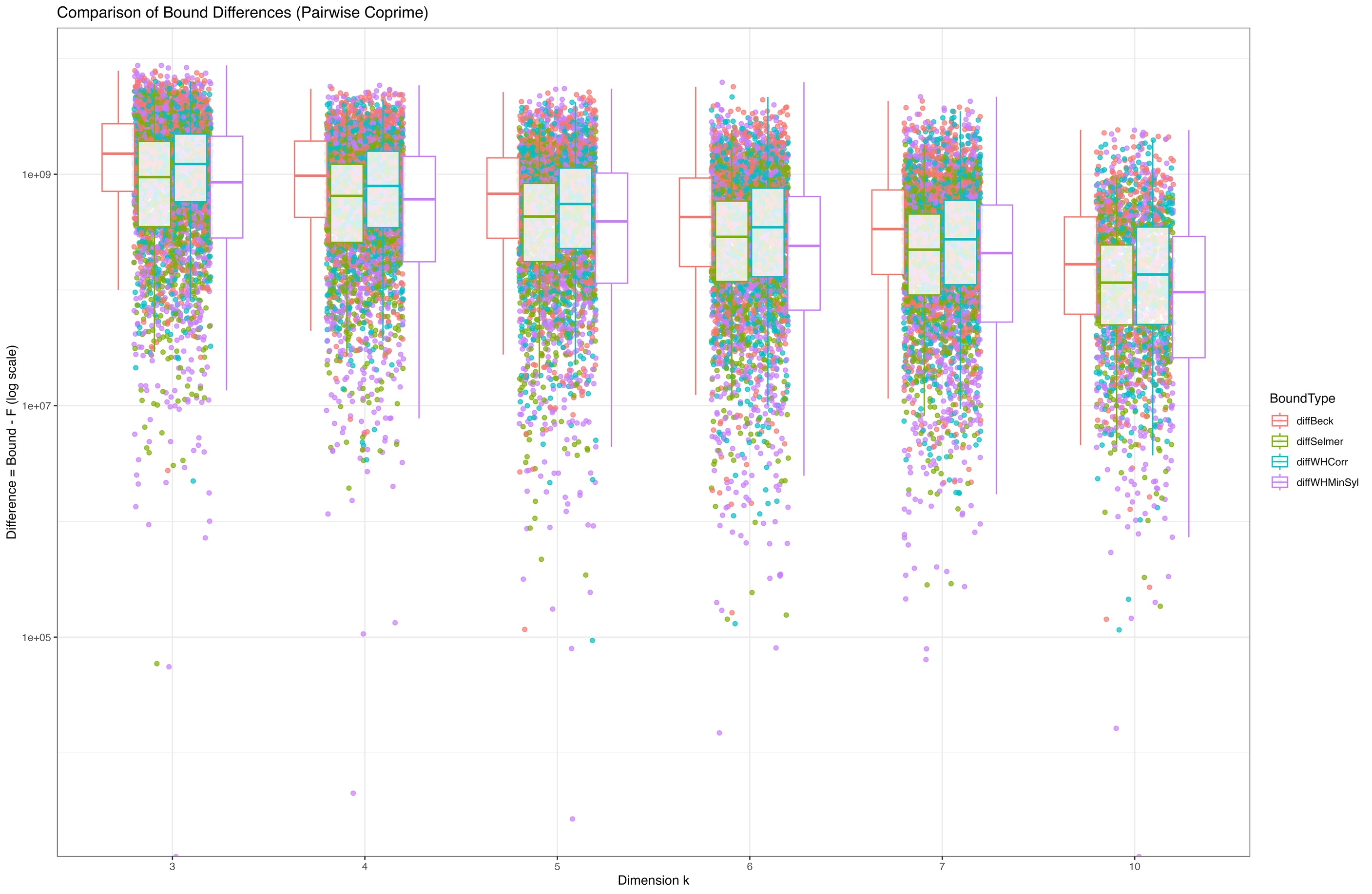}
\caption{Box plots of the difference, i.e. $\text{Bound} - F(\boldsymbol{a})$, across dimensions with $\| \boldsymbol{a} \|_{\infty} \le 100,000$. The log-scaled vertical axis reveals the comparative tightness and spread of different bounds.}
\label{fig:boxplot_all_k}
\end{figure}

\vspace{2.0mm}

Figure~\ref{fig:boxplot_all_k} displays log-scaled box plots of the differences, i.e. $\text{Bound} - F(\boldsymbol{a})$, for each bound type, grouped by dimension. Notice that error magnitudes generally decrease as the dimension $n$ increases, which aligns with the expectation of tighter asymptotics in higher dimensions. The second bound of Williams and Haijima \eqref{williams_sylvester_extension} consistently achieves the lowest median errors across dimensions. Moreover, its spread tends to extend downward compared to other bounds, highlighting its somewhat frequent tendency to provide tighter approximations. The bound of Selmer \eqref{Selmer_bound_n_1}, in contrast, demonstrates tighter and more symmetric distributions, suggesting more robust performance and predictable performance across varying dimensions. The bound of Beck et al. \eqref{Beck_et_al_bound_1}, however, exhibits notably wider dispersion and consistently higher median errors, reflecting its typically weaker tightness. A corresponding box plot for input vectors with $\|\boldsymbol{a}\|_\infty \le 10{,}000$ is provided in Appendix~\ref{appendix:box_m10000}, which demonstrates similar trends on a smaller input scale.

\vspace{2.0mm}

\begin{figure}[ht!]
\centering
\includegraphics[width=1\linewidth]{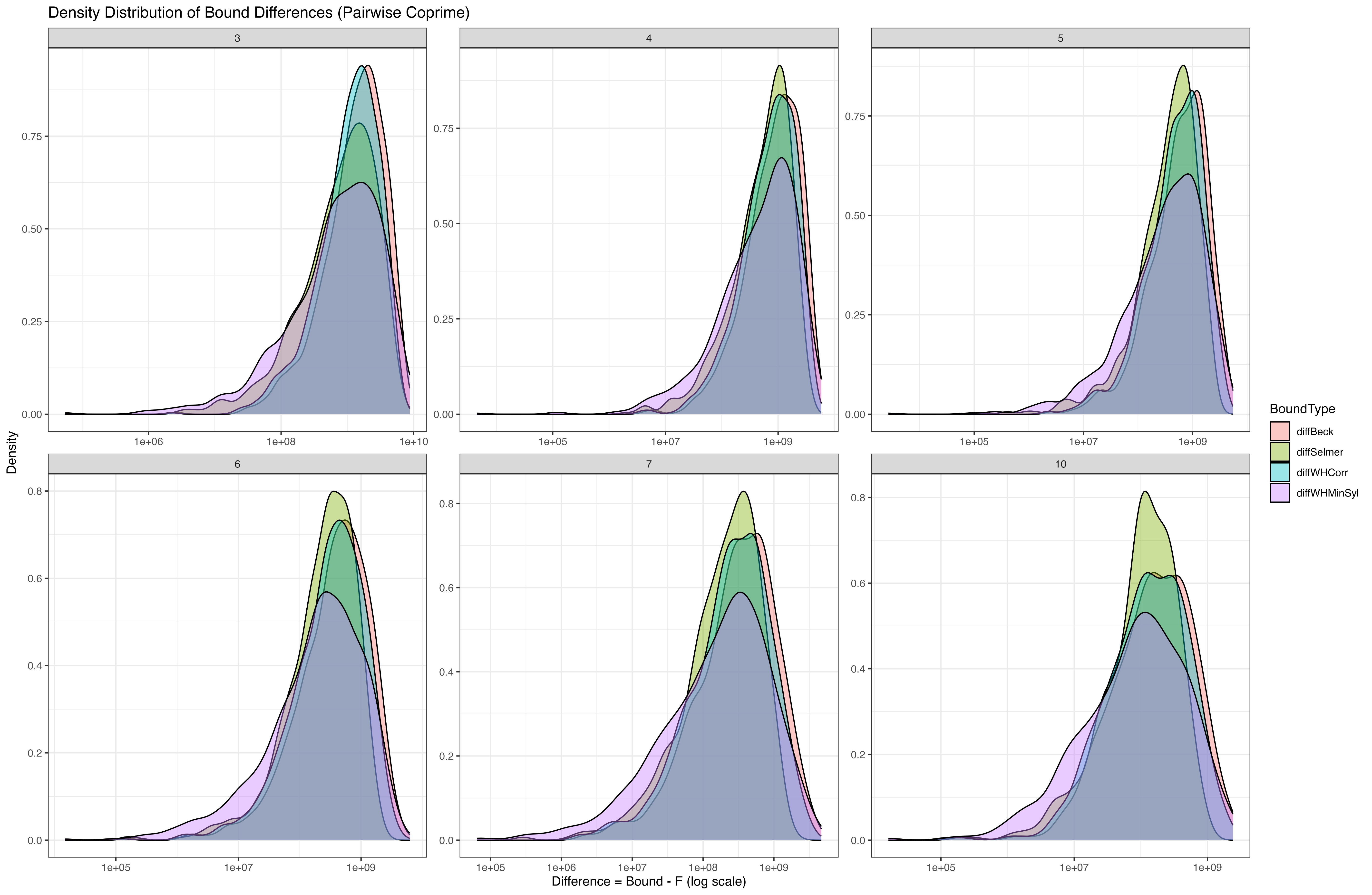}
\caption{Density plots of $\text{Bound} - F(\boldsymbol{a})$ for $\| \boldsymbol{a} \|_{\infty} \le 100,000$ across dimensions $n$. The log-scaled vertical axis reveals the comparative tightness and spread of different bounds.}
\label{fig:density_by_k}
\end{figure}

While box plots succinctly summarise distribution characteristics, they obscure the precise shape of the underlying distribution. To investigate this, we use kernel density estimates (KDEs). Note for completeness that KDEs provides a smooth approximation of the underlying probability density function by aggregating \enquote{bumps} (kernels) placed at each data point. This results in a continuous curve that highlights where data points tend to cluster. In contrast to classical density models that assume a specific shape (e.g. Gaussian), KDEs offer a flexible and data-driven way to visualise the distribution without strong parametric assumptions. Figure~\ref{fig:density_by_k} presents KDE plots for bound differences. 
It should be noted that each KDE curve peaks in regions where the observed differences are most densely concentrated, where the areas under curves sum to one. 
We observe that the second bound of Williams and Haijima \eqref{williams_sylvester_extension} (purple) demonstrates a slight leftward shift relative to other bounds, indicating that its most frequent errors are smaller. Its distribution does however have a slightly longer right tail compared to some bounds, highlighting occasional large overestimations despite its typically sharper predictions. The bound of Selmer \eqref{Selmer_bound_n_1} (green) displays a more pronounced peak, reflecting consistent accuracy and less frequent large deviations. The bound of Beck et al. \eqref{Beck_et_al_bound_1} (orange) has a somewhat wide distribution that is slightly rightward shifted, signifying generally weaker performance. This outcome aligns with theoretical expectations given that \eqref{Beck_et_al_bound_1} is known analytically to be tighter than \eqref{Williams_corrected_1} only in a (small) finite number of cases \cite[Theorem 3]{williams2023considering}. For a similar comparison with the smaller upper bound $\|\boldsymbol{a}\|_\infty \le 10,000$, kernel density plots are provided in Appendix~\ref{appendix:density_structured_m10000}.

\vspace{2.0mm}

\vspace{2.0mm}

\begin{figure}[ht!]
\centering
\includegraphics[width=1\linewidth]{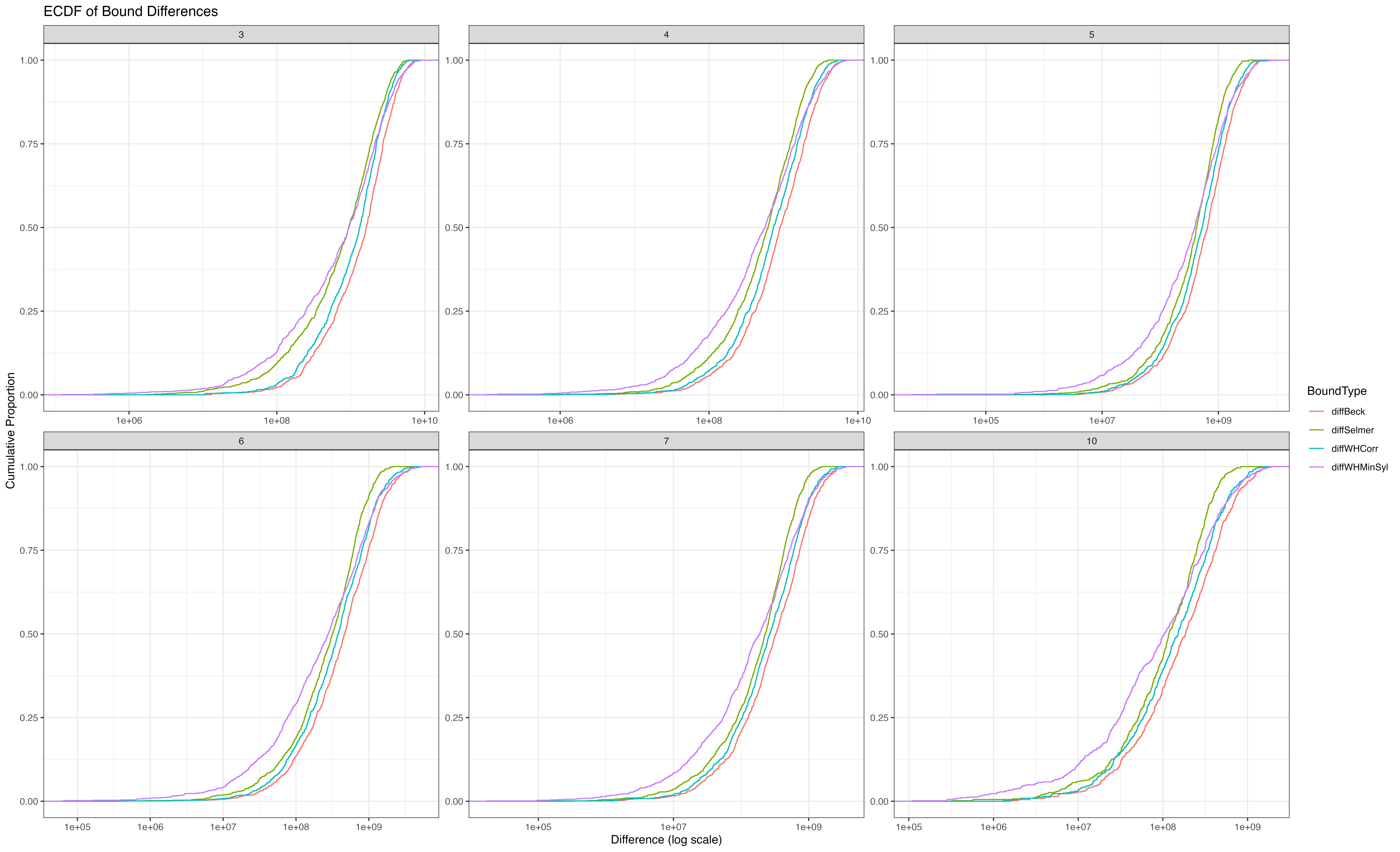}
\caption{Empirical cumulative distribution functions of bound differences across dimensions $n$ for $\| \boldsymbol{a} \|_{\infty} \le 100,000$.}
\label{fig:ecdf_by_k}
\end{figure}

While KDEs reveal underlying distribution shapes, the empirical cumulative distribution functions (ECDFs) in Figure~\ref{fig:ecdf_by_k} provide a complementary view by illustrating the proportion of observations below a given error threshold. Each curve increases (in the vertical direction) monotonically from zero to one, capturing the cumulative performance of each bound across all quantiles.
Note that interpreting an ECDF is straightforward, namely for a given error value on the horizontal axis, a higher curve indicates better performance since, in such case, it corresponds to a greater proportion of instances with errors less than or equal to that value. 
The ECDFs in Figure~\ref{fig:ecdf_by_k} confirm the strong performance of the second bound of Williams and Haijima \eqref{williams_sylvester_extension} across most quantiles. The bound of Selmer \eqref{Selmer_bound_n_1} ranks second in performance over the majority of the distribution, while the other bounds perform less well across the majority of their range, as reflected by their lower ECDF curves.

\vspace{2.0mm}

\vspace{2.0mm}

\begin{figure}[ht!]
\centering
\includegraphics[width=1\linewidth]{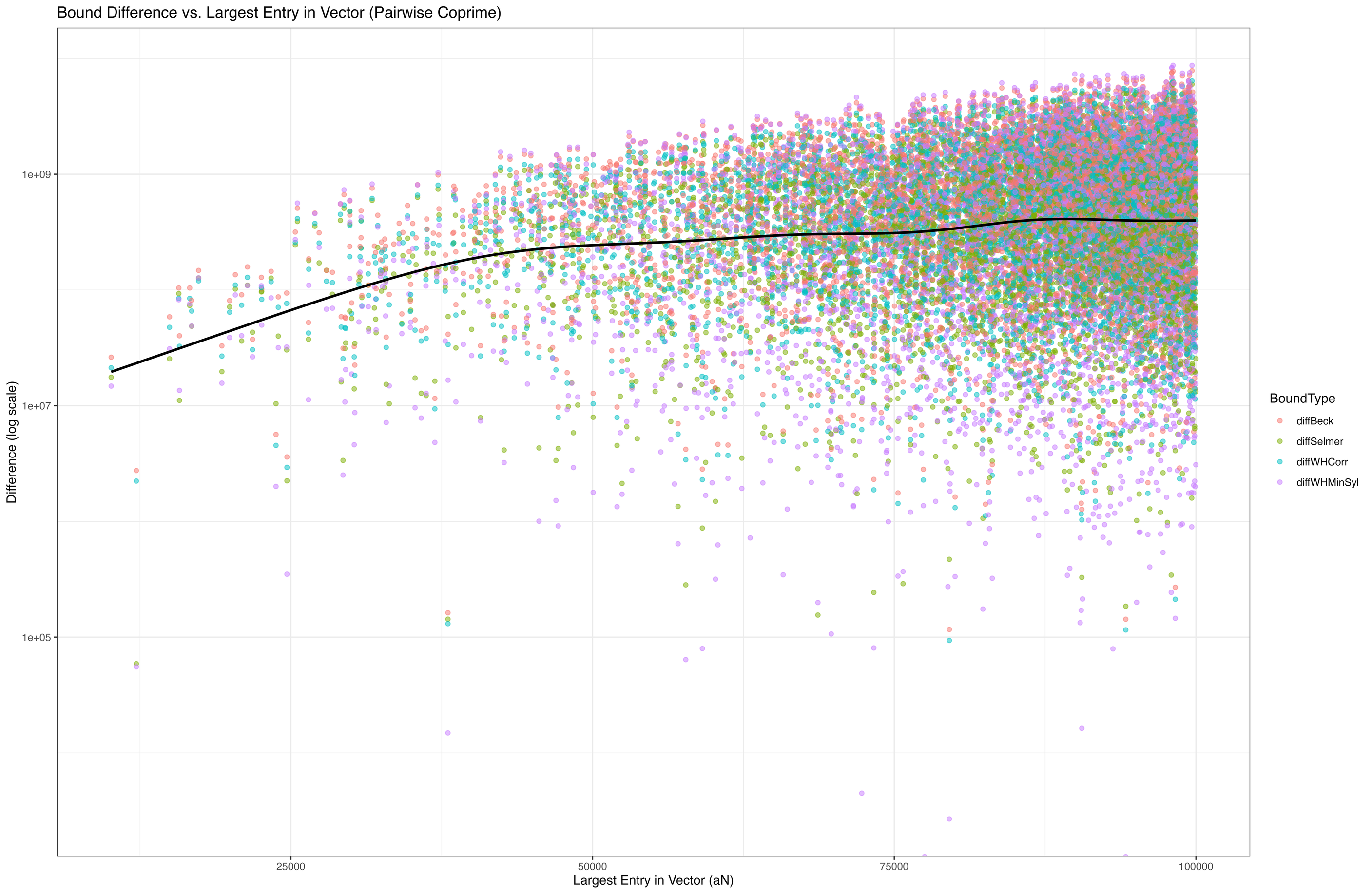}
\caption{Scatter plot of the difference, i.e. $\text{Bound} - F(\boldsymbol{a})$, against the maximum entry $a_n$ separated by bound type, with GAM smoother for $\| \boldsymbol{a} \|_{\infty} \le 100,000$.}
\label{fig:scatter_diff_vs_an}
\end{figure}

\begin{figure}[ht!]
\centering
\includegraphics[width=1\linewidth]{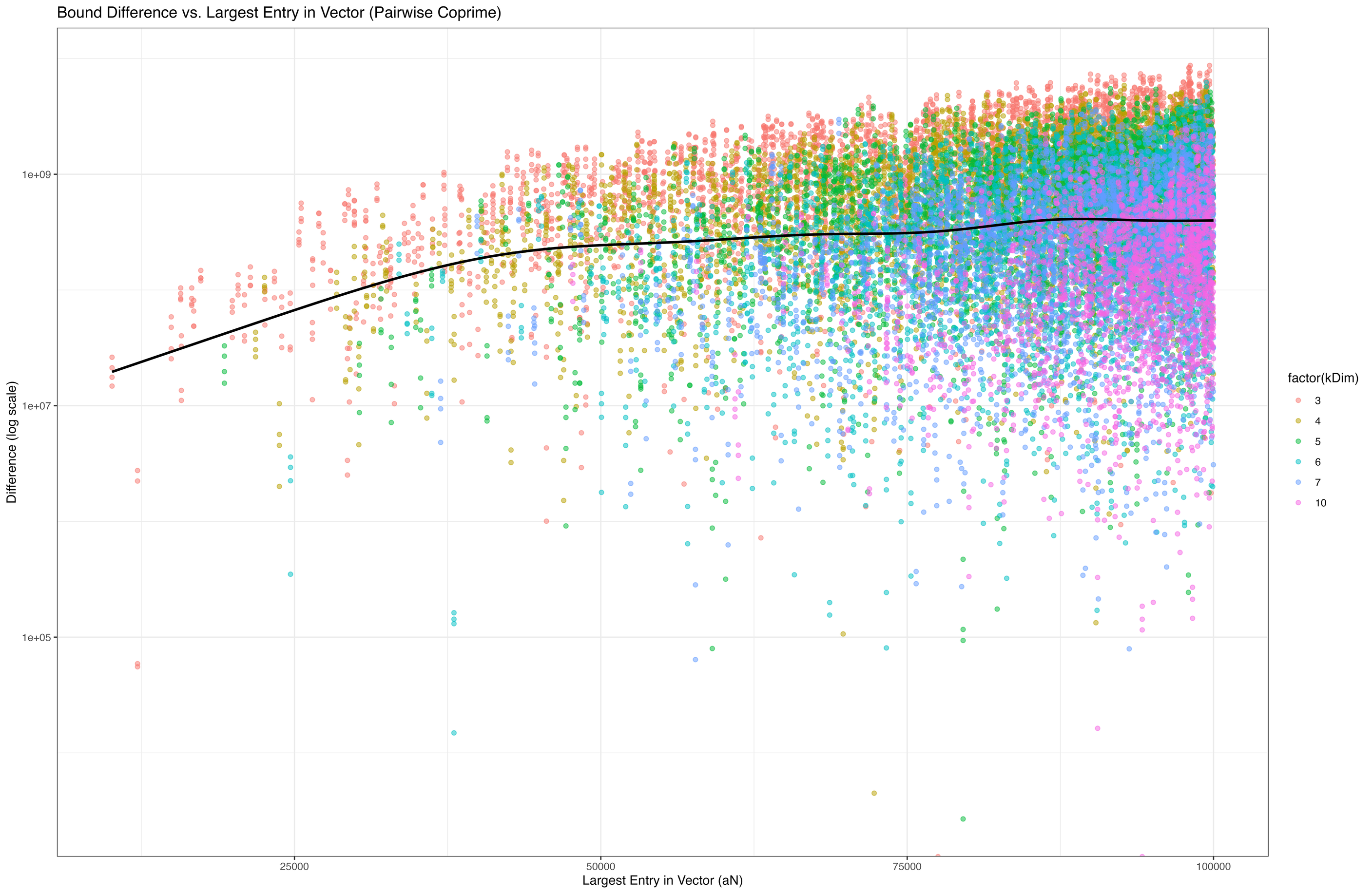}
\caption{Scatter plot of the difference, i.e. $\text{Bound} - F(\boldsymbol{a})$, against the maximum entry $a_n$ separated by dimension, with GAM smoother.}
\label{fig:scatter_diff_vs_an_dim}
\end{figure}

Figures \ref{fig:scatter_diff_vs_an} and \ref{fig:scatter_diff_vs_an_dim} show scatter plots of the bound differences against the maximum absolute valued entry $\| \boldsymbol{a} \|_{\infty} = a_n$. It should be noted that each point represents a simulation instance, coloured by either the bound type (Figure \ref{fig:scatter_diff_vs_an}) or by the dimension (Figure \ref{fig:scatter_diff_vs_an_dim}). The black curve is a smoothed generalized additive model (GAM) curve, which highlights the general scaling trend against increases in $a_n$. These figures suggest that bound errors tend to increase as the largest vector entry $a_n$ increases, although the rate of increase appears to slow. There remains considerable variability around the trend, especially for large values of $a_n$. 

\vspace{2.0mm}

\vspace{2.0mm}

\vspace{2.0mm}

Figure \ref{fig:scatter_diff_vs_an} clarifies differences among the various bounds. The bound of Beck et al. \eqref{Beck_et_al_bound_1} consistently exhibits large errors across the range of $a_n$, aligning with previous observations. The bounds of Selmer \eqref{Selmer_bound_n_1} and the first bound of Williams and Haijima \eqref{Williams_corrected_1} appear relatively similar across all $a_n$. The second bound of Williams and Haijima \eqref{williams_sylvester_extension} frequently achieves the lowest errors, particularly for moderate to large $a_n$, as evidenced by its concentration in the lower region of the vertical axis. 
Despite this, it should be emphasised that this bound also exhibits substantial dispersion, attaining the highest observed errors among the bounds in a non-negligible number of random instances. This behaviour underscores a trade-off, namely that while the bound \eqref{williams_sylvester_extension} tends to yield the best performance in some cases, its wide variability renders it less reliable in worst-case settings. A corresponding plot for input vectors with $\| \boldsymbol{a} \|_{\infty} \le 10,000$ is provided in Appendix \ref{appendix:boxplot_m10000_structured}.

\vspace{2.0mm}

Figure~\ref{fig:scatter_diff_vs_an_dim} shows the same bound error differences plotted against $a_n$, but where the colours now indicate the dimension $n$, rather than the bounding function. Observe that for any fixed $a_n$, higher dimensional vectors (e.g. $n=10$, shown in pink) seem to produce smaller absolute errors. This is consistent with the theoretical expectation that many bounds improve with the dimension due to asymptotic behaviour on average. Despite this, the large overlap between dimensions and the considerable vertical spread of points reinforces that the dimension alone does not fully explain error magnitude. 

\vspace{2.0mm}

\begin{figure}[ht!]
\centering
\includegraphics[width=1\linewidth]{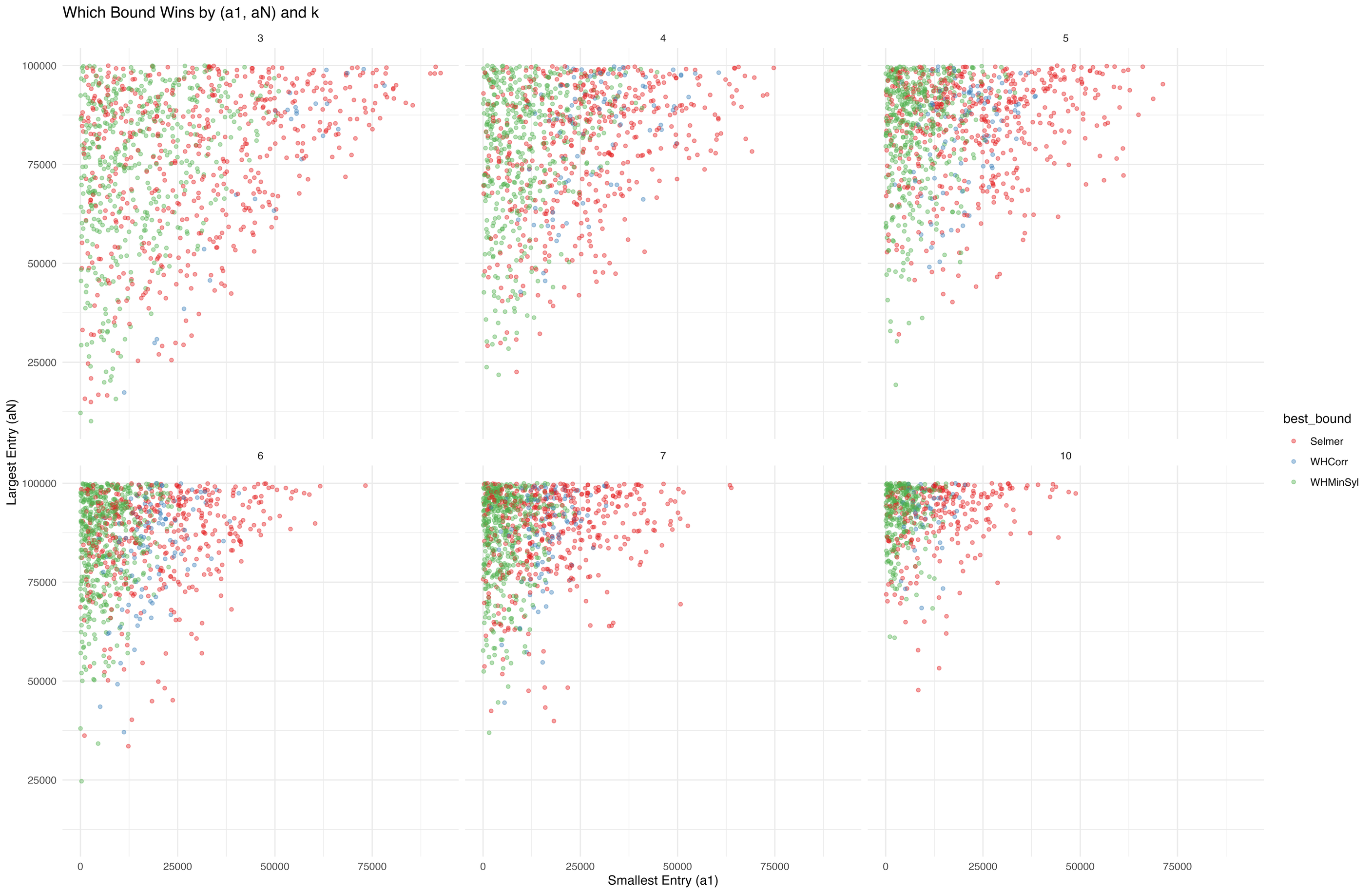}
\caption{Scatter plot of best-performing bound across $(a_1, a_n)$.}
\label{fig:scatter_best_bound_by_inputs}
\end{figure}

We now identify which bound is tightest across different regions of the input space. In Figure~\ref{fig:scatter_best_bound_by_inputs} each point corresponds to one simulated instance, plotted by its smallest and largest entries $(a_1, a_n)$, and coloured by the bound that minimises the error. Note that for clarity we omit the bound of Beck et al. \eqref{Beck_et_al_bound_1} here since it is known to be tighter than \eqref{Williams_corrected_1} in only a finite number of cases \cite[Theorem 3]{williams2023considering}. The second bound of Williams and Haijima \eqref{williams_sylvester_extension} (green) dominates the lower triangle of the input space, corresponding to small $a_1$ and moderate $a_n$. The bound of Selmer \eqref{Selmer_bound_n_1} (red) dominates in the upper right of the input space, i.e. when $a_n$ is large and $a_1$ is moderate. The first bound of Williams and Haijima \eqref{Williams_corrected_1} (blue) maintains a consistent presence across intermediate regions. These clusters reflect how structural features of the input, such as scale and conditioning, influence which bound is tightest. An analogous regional performance plot for smaller inputs can be found in Appendix~\ref{appendix:best_bound_region_m10000}. Further, we present a summary of the above outcomes, showing the overall proportion of instances for which each bound achieves the tightest estimate, in Appendix~\ref{appendix:best_bound_proportions}.

\vspace{2.0mm}

\vspace{2.0mm}

To sharpen the focus of our analysis, we now restrict attention to the bounds of Selmer \eqref{Selmer_bound_n_1} and the first bound of Williams and Haijima \eqref{Williams_corrected_1}. Note that while the second bound of Williams and Haijima \eqref{williams_sylvester_extension} frequently provides the tightest approximation in typical cases, as seen in Figures \ref{fig:density_by_k} and \ref{fig:ecdf_by_k}, it exhibits substantial variability and can produce larger errors in extreme cases. The bound of Beck et al. \eqref{Beck_et_al_bound_1}, on the other hand, is known to be tighter than  \eqref{Williams_corrected_1} in only a finite number of cases. Thus, it is natural to compare the bound of Selmer \eqref{Selmer_bound_n_1} and the first bound of Williams and Haijima \eqref{Williams_corrected_1} given their stable error profiles, moderate variability, and analytical foundations.

\vspace{2.0mm}

\begin{figure}[ht!]
\centering
\includegraphics[width=1\linewidth]{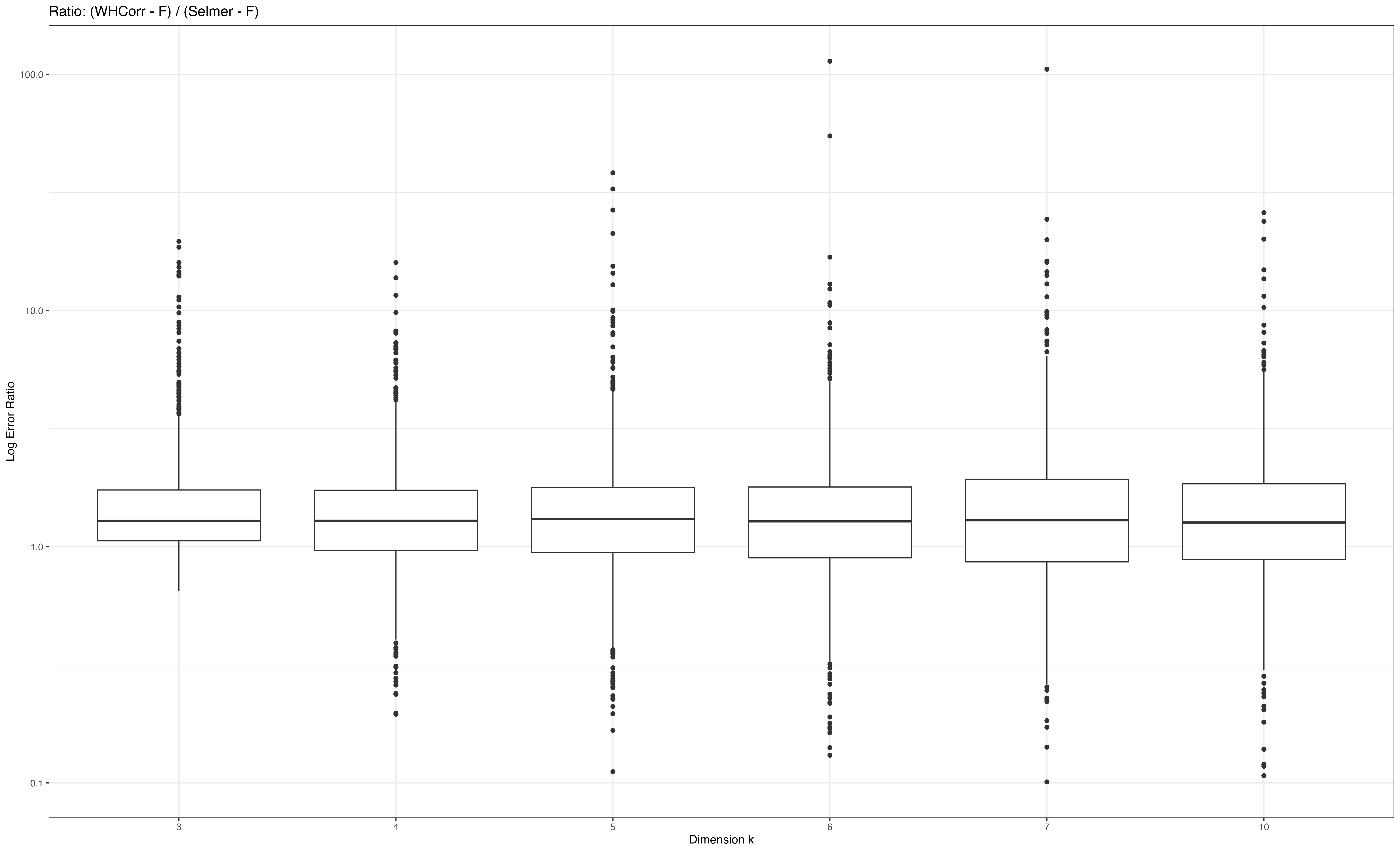}
\caption{Log-scaled box plots of the ratio between the first bound of Williams and Haijima \eqref{Williams_corrected_1} and the bound of Selmer \eqref{Selmer_bound_n_1}, across dimensions. Note that a ratio greater than 1 indicates \eqref{Selmer_bound_n_1} is tighter, while a ratio less than 1 indicates that \eqref{Williams_corrected_1} is tighter.}
\label{fig:WHCorr_vs_Selmer_ratio}
\end{figure}

\vspace{2.0mm}

Figure \ref{fig:WHCorr_vs_Selmer_ratio} presents log-scaled box plots of the ratio of absolute errors between the (corrected) Williams and Haijima bound \eqref{Williams_corrected_1} and the bound of Selmer \eqref{Selmer_bound_n_1} across dimensions. In particular, the plotted quantity is 
$$
R_n = \frac{\texttt{diffWHCorr} -F(\boldsymbol{a})}{\texttt{diffSelmer} -F(\boldsymbol{a})}
$$
across dimensions $n$. Note that a ratio $R_n > 1$ means that the bound of Selmer \eqref{Selmer_bound_n_1} has the smaller absolute error, while $R_n < 1$ means that the bound of Williams and Haijima bound \eqref{Williams_corrected_1} is tighter. The medians lie just above 1 across all dimensions, suggesting that the bound of Selmer \eqref{Selmer_bound_n_1} typically provides a slightly tighter bound on $F(\boldsymbol{a})$. Despite this, the bulk of the distribution is concentrated around the median, indicating that \eqref{Williams_corrected_1} performs competitively in most cases.
A small number of significant outliers appear in each dimension, highlighting instances where \eqref{Williams_corrected_1} considerably overestimates $F(\boldsymbol{a})$. It is interesting that in higher dimensions we observe a growing number of outliers in the opposite direction, where \eqref{Selmer_bound_n_1} significantly \enquote{overshoots} relative to \eqref{Williams_corrected_1}.

\vspace{2.0mm}

\begin{figure}[ht!]
\centering
\includegraphics[width=1\linewidth]{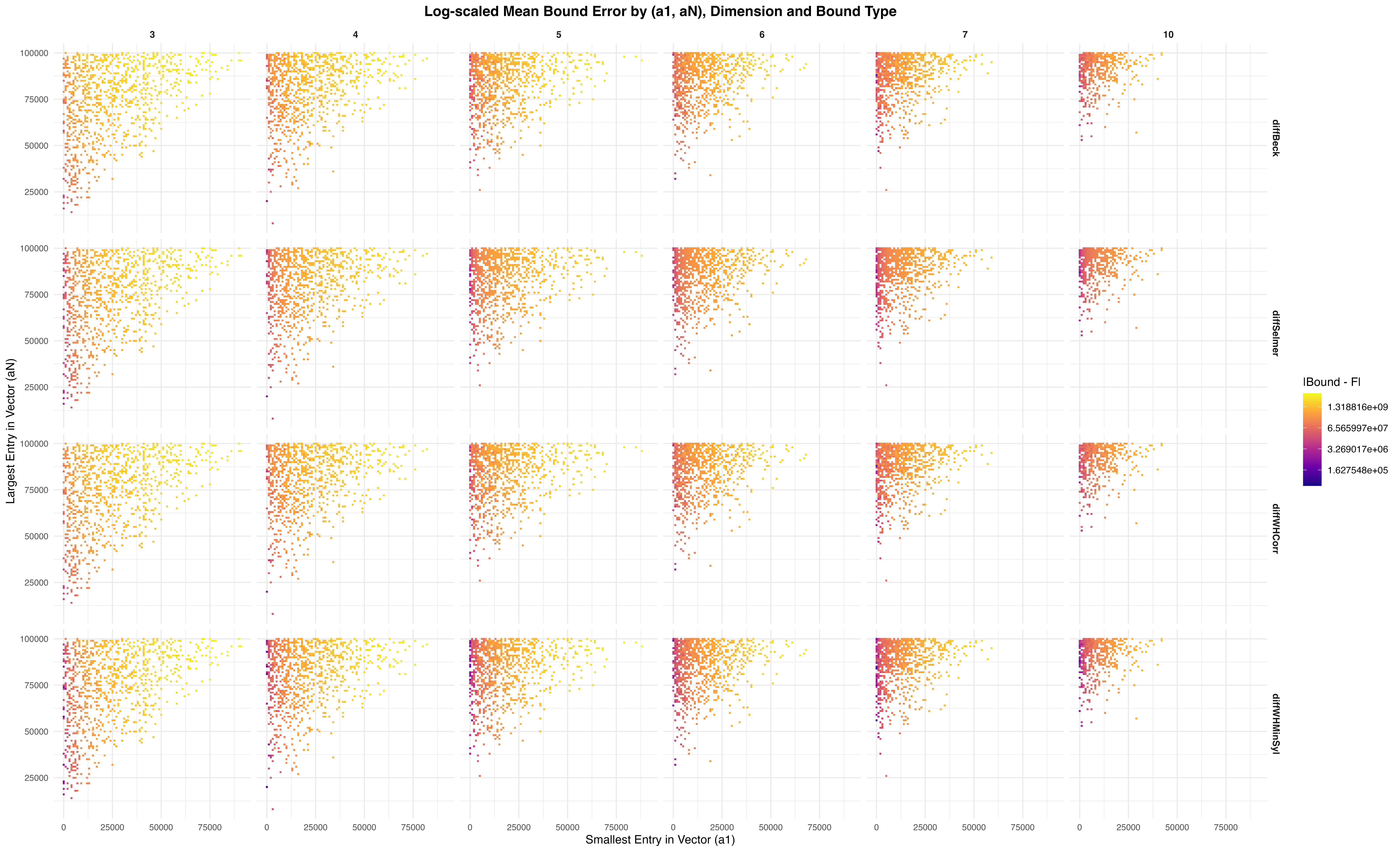}
\caption{Bound errors plotted by smallest and largest vector entries $(a_1, a_n)$. Light colours (yellow) indicate larger overestimates, while darker colours (purple) indicates smaller errors.}
\label{fig:heatmap_all_bounds_structured}
\end{figure}

To better understand how the structural characteristics of input vectors impact upon the tightness of different bounds, we illustrate the absolute error across different combinations of the smallest and largest entries, namely $a_1$ and $a_n$. Figure~\ref{fig:heatmap_all_bounds_structured} plots each randomly sampled vector $\boldsymbol{a}$ according to its smallest and largest entries $(a_1, a_n)$, coloured by the corresponding absolute error of each bound. The colour scale is log-transformed to effectively visualise the wide range of observed error magnitudes.
The plots confirm that increasing the largest entry $a_n$ tends to correspond to larger bound errors, consistent with expected asymptotic behaviour. 
We also observe a consistent decrease in errors which increasing dimension (evidenced by the darkening colour profile), a trend that is consistent across bounds. The first bound of Williams and Haijima \eqref{Williams_corrected_1} tends to achieve the lowest errors when $a_1$ is small, regardless of $a_n$. The bound of Selmer \eqref{Selmer_bound_n_1} is likewise most accurate for small values of $a_1$, although its error profile degrades slightly more gradually across the input domain. In contrast, \eqref{williams_sylvester_extension} exhibits high volatility, with good performance for small $a_1$, but sharp increases in error magnitude as $a_1$ becomes large.
Appendix~\ref{appendix:mean_error_heatmap_m10000} presents the same heat map visualisation for inputs with $\|\boldsymbol{a}\|_\infty \le 10,000$, which illustrates similar structural patterns with slightly reduced error magnitudes.

\vspace{2.0mm}

\vspace{2.0mm}

\begin{figure}[ht!]
\centering
\includegraphics[width=0.9\linewidth]{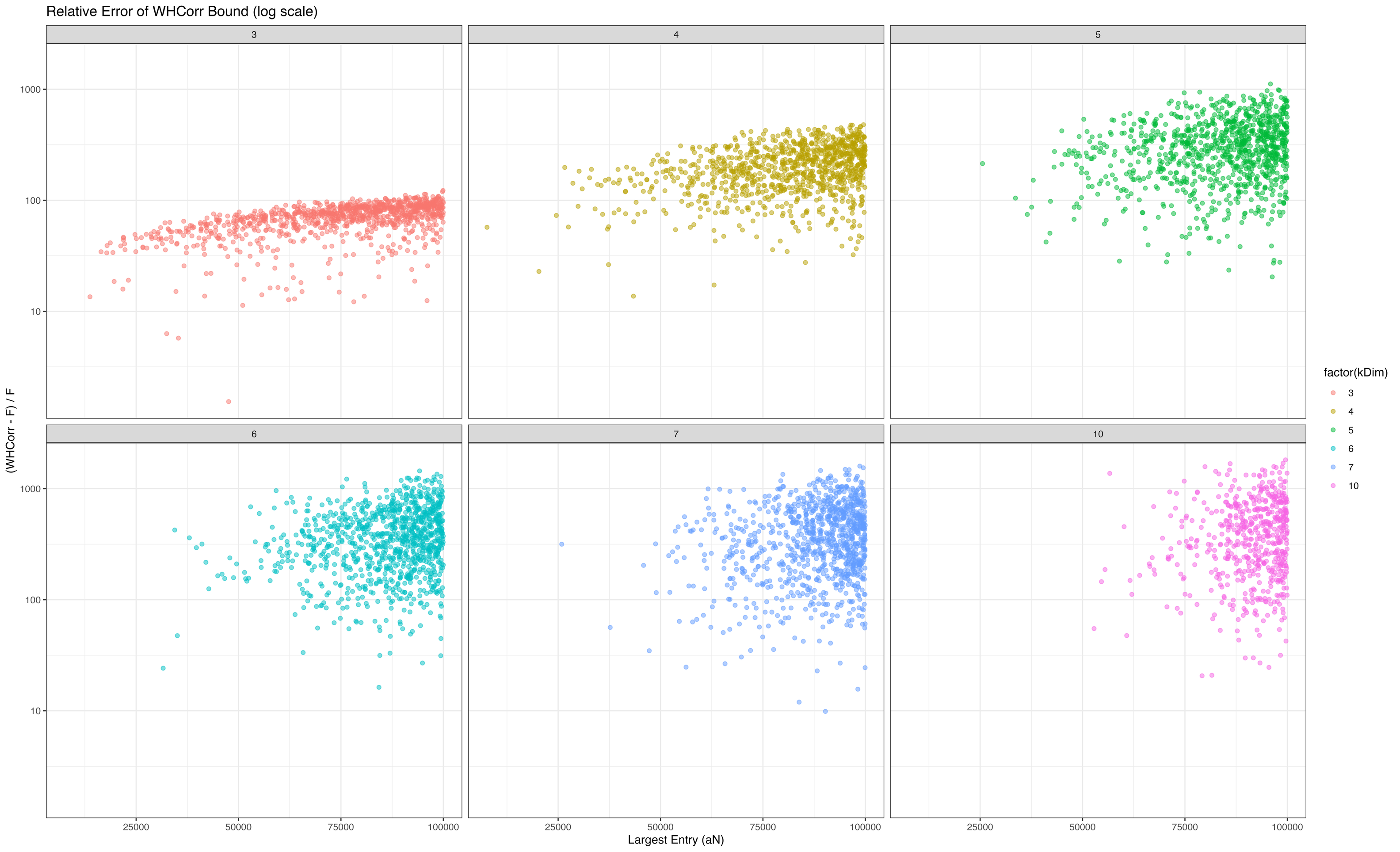}
\caption{Relative error $\frac{\texttt{WHCorr} - F}{F}$ plotted against $a_n$ across dimensions.}
\label{fig:rel_err_whcorr_vs_an}
\end{figure}

To investigate how the first bound of Williams and Haijima \eqref{Williams_corrected_1} scales with input size, we examine its relative error, i.e. the amount by which the bound exceeds the Frobenius number expressed as a proportion of that number. Figure~\ref{fig:rel_err_whcorr_vs_an} plots this quantity against the maximum entry $a_n$, grouped by dimension $n$. The vertical axis is presented on a logarithmic scale to accommodate both moderate and extreme deviations across the input space. The plots reveals that in lower dimensions, the bound \eqref{Williams_corrected_1} displays relatively modest and slowly increasing relative errors as the maximum entry $a_n$ grows. In contrast, for higher dimensional settings, the relative error exhibits a more pronounced upward trajectory. This suggests that while the bound \eqref{Williams_corrected_1} performs robustly for small-scale instances, its proportional overestimation deteriorates more noticeably in higher dimensions. For a more detailed analysis of how input imbalance influences bound quality, Appendix~\ref{appendix:WH_Error_Input} examines how the relative error in \eqref{Williams_corrected_1} varies with the    ratio $a_n / a_1$.

\vspace{2.0mm}

\section{Limitations of Upper Bounds}
Recall that we have so far compared the relative tightness of known upper bounds on the Frobenius number $F(\boldsymbol{a})$ using Monte Carlo simulation techniques. Notably, the bounds that are valid under the (weaker) conditions \eqref{conditions on a} share the characteristic of being quadratic in the worst-case and depend on $\| \boldsymbol{a} \|_{\infty}$. Indeed, Williams and Haijima \cite[Lemma 1]{williams2023considering} prove that any general upper bound on $F(\boldsymbol{a})$ must inherently depend on the maximum absolute entry $\| \boldsymbol{a} \|_{\infty}$. A natural question then arises, namely is it possible to identify a general upper bound for $F(\boldsymbol{a})$  under the (weaker) conditions \eqref{conditions on a}, but with a worst-case exponent lower than quadratic? The following result formally demonstrates that this is not feasible, where the proof is provided thereafter.

\vspace{2.0mm}

\begin{theorem} \label{quadratic_infeasiblity}
Fix any integer $n \ge 2$. Suppose an integral vector $\boldsymbol{a}=\left(a_1, a_2, \ldots, a_n\right)^T$ satisfies \eqref{conditions on a} with $a_1 \le a_2 \le \cdots \le a_n$. Then, for every $\varepsilon>0$, there is no constant $C > 0$ (independent of $\boldsymbol{a}$) such that 
$$
F(\boldsymbol{a}) \le\ C \cdot \big( a_{1}\,a_{n} \big)^{1-\varepsilon}
$$
holds for all vectors $\boldsymbol{a}$ satisfying \eqref{conditions on a}.
\end{theorem}

\vspace{2.0mm}

It should be emphasised that the above indicates that one cannot achieve an asymptotic upper bound with order less than $a_1 a_n$ unless additional conditions are imposed. In other words, this result demonstrates that one cannot globally force the Frobenius number $F(\boldsymbol{a})$ to be bounded \enquote{universally} by a sub-linear power of $a_1 a_n$ without imposing further restrictions.

\vspace{2.0mm}


\section{Proof of Proposition \ref{prop1}}
\begin{proof} 
Observe that if $n=2$, then the stronger \eqref{stronger conditions on a} and weaker conditions \eqref{conditions on a} are equivalent. Thus, we focus only on the case $n \ge 3$, where we show that there are counterexamples for any dimension $n$.

\vspace{2.0mm}

Let us consider two cases, namely $n=3$ and $n\ge4$, respectively. If $n=3$, then consider the vector $\boldsymbol{a} = (8,32,59)^T$ with $F(\boldsymbol{a}) = 405$. In this case, the bound of Selmer \eqref{Selmer_bound_n} yields
$$
2 \cdot 59 \Big\lfloor{\,\frac{8}{3}\,\Big\rfloor} - 8 = 2 \cdot 59 \cdot 2 - 8 = 228
$$
and clearly $405 \not\le 228$. Thus, the upper bound \eqref{Selmer_bound_n} fails when $n=3$. 

\vspace{2.0mm}

If, instead, $n\ge4$, then let
$$
\widetilde{\boldsymbol{a}} = \big(8, \, 8 k_2, \, 8k_3, \, \ldots, \, 8k_{n-1}, \, 59\big)^T
$$ 
satisfy $k_i \in \mathbb{Z}_{>0}$ for each $i$ with $k_2 \le k_3 \le \cdots \le k_{n-1} \le \lfloor 59/8 \rfloor =  7$ provides a counterexample for any $n$ as required. Indeed, by construction, we have $F(\boldsymbol{a}) = F(\widetilde{\boldsymbol{a}}) = 405$, whereas the upper bound \eqref{Selmer_bound_n} 
becomes
$$
2 \cdot 59 \Big\lfloor{\,\frac{8}{n}\,\Big\rfloor} - 8 \le 2 \cdot 59 \Big\lfloor{\,\frac{8}{3}\,\Big\rfloor} - 8 = 228,
$$
where the first inequality follows since $\lfloor{8/n}\rfloor \le \lfloor{8/3}\rfloor$ holds for $n \ge 3$. Consequently, clearly $405 \not\le 228$, and the bound fails for any $n \ge 4$, which completes the proof. 
\end{proof}

\vspace{2.0mm}

It should be emphasised that the choice of $\boldsymbol{a}=(8,32,59)^T$ in the above is not unique. In particular, Table \ref{table:comparison} presents other vectors for which the bound \eqref{Selmer_bound_n} fails for $n=3$ under the weaker conditions \eqref{conditions on a}, where a similar embedding argument could be presented.

\vspace{2.0mm}

\begin{table}[ht!]
    \centering
    \begin{subtable}{0.48\textwidth}
        \centering
        \begin{tabular}{c | c | c}
        $\boldsymbol{a} = (a_1, a_2, a_3)^T$ & $F(\boldsymbol{a})$ & Bound (Selmer) \\
        \hline
        $(4, 12, 25)^T$  & 71  &   46 \\
        $(4, 24, 31)^T$ &  89 &  58 \\
        $(4, 32, 57)^T$ & 167 &  110 \\
        $(4, 39, 52)^T$ &  113 &  100 \\
        $(4, 43, 44)^T$ &  125 &  84 \\
        $(4, 44, 45)^T$ &  131  &   86 \\
        \hline
        $(5,  7, 12)^T$ &  23  &   19 \\
        $(5, 10, 33)^T$ & 127 &  61 \\
        $(5, 13, 20)^T$ &  47   &   35 \\
        $(5, 15, 31)^T$ & 119   &   57 \\
        $(5, 16, 20)^T$ &  59    &  35 \\
        $(5, 24, 34)^T$  & 91    &  63 \\
        $(5, 28, 50)^T$ &107   &   95 \\
        $(5, 30, 39)^T$ & 151    &  73 \\
        $(5, 30, 41)^T$ & 159   &   77 \\
        $(5, 31, 50)^T$   &119    &  95 \\
        $(5, 32, 37)^T$ &  123   &   69 \\
        $(5, 34, 39)^T$ &  131    &  73 \\
        $(5, 37, 40)^T$ &  143    &  75 \\
        $(5, 38, 58)^T$ &  147   &  111 \\
        $(5, 45, 53)^T$ &  207   &  101 \\
        $(5, 46, 51 )^T$ & 179   &   97 \\
        $(5, 47, 55)^T$ &  183   &  105 \\
        $(5, 48, 58)^T$ &  187  &   111 \\
        \end{tabular}
        \caption{Cases with $a_1 \in \{4,5\}$}
        \label{table_of_examples_a1_45}
    \end{subtable}
    \hfill
    \begin{subtable}{0.48\textwidth}
        \centering
        \begin{tabular}{c | c | c}
        $\boldsymbol{a} = (a_1, a_2, a_3)^T$ & $F(\boldsymbol{a})$ & Bound (Selmer) \\
        \hline
        $(6, 18, 47)^T$ &  229   &  182 \\
        $(6, 29, 30)^T$ &  139  &   114 \\
        \hline
        $(7 ,14, 44)^T$ &  257  &   169 \\
        $(7 ,21, 22)^T$ &  125  &    81 \\
        $(7 ,21, 50)^T$ &  293  &   193 \\
        $(7 ,37, 44)^T$ &  215   &  169 \\
        $(7 ,42, 57)^T$ &  335   &  221 \\
        $(7 ,42, 60 )^T$ & 353  &   233 \\
        $(7 ,48, 55)^T$ &  281  &   213 \\
        $(7 ,51, 58 )^T$ & 299   &  225 \\
        \hline
        $(8, 16, 55)^T$ &  377 &   212 \\
        $(8, 23, 40)^T$ &  153  & 152 \\
        $(8, 24 ,31)^T$ &  209 &   116 \\
        $(8, 24 ,41)^T$ &  279 &    156 \\
        $(8, 24 ,49)^T$ &  335  &   188 \\
        $(8, 31 ,40)^T$ &  209  &   152 \\
        $(8, 32, 59)^T$ &  405  &   228 \\
        $(8, 39 , 48)^T$ &  265  &  184
        \end{tabular}
        \caption{Cases with $a_1 \in \{6,7,8\}$}
        \label{table_of_examples_a1_678}
    \end{subtable}
\caption{Comparison of the Frobenius number \( F(\boldsymbol{a}) \) and the Selmer \cite{selmer1977linear} bound \eqref{Selmer_bound_n} for various integer triples $\boldsymbol{a}=(a_1, a_2, a_3)^T$ satisfying the (weaker) conditions \eqref{conditions on a}, which are for convenience constrained by $\| \boldsymbol{a} \|_{\infty} = \max_i |a_i| \le 60$ and \( a_1 \in \{4,5,6,7,8\} \). The left table focuses on cases where $a_1 \in \{4,5\}$, while the right table includes cases where $a_1 \in \{6,7,8\}$.}
\label{table:comparison}
\end{table}


\section{Proof of Theorem \ref{quadratic_infeasiblity}}
Observe that it is sufficient to assume that $\varepsilon<1$ holds. In particular, if $\varepsilon=1$, then we yield
$$
\bigl(a_{1}a_{n}\bigr)^{\,1-\varepsilon} =
\bigl(a_{1}a_{n}\bigr)^{0} = 1.
$$
Thus, the claimed bound in Theorem \ref{quadratic_infeasiblity} becomes 
$$
F(\boldsymbol{a}) \le\ C
$$
for all $\boldsymbol{a}$ satisfying \eqref{conditions on a}. However, the Frobenius number grows arbitrarily large as $\min(a_i) \rightarrow \infty$ (see e.g. \cite[Proposition 2]{williams2023considering}), which eventually exceeds any constant $C$. 
If instead $\varepsilon >1$, then $1-\varepsilon < 0$ holds, and the claimed bound becomes
\begin{equation} \label{varepsilon<0_case}
\begin{aligned}
F(\boldsymbol{a}) \le C \cdot \big( a_{1}\,a_{n} \big)^{1-\varepsilon} \\
                    = \frac{C}{\big( a_{1}\,a_{n} \big)^{\varepsilon-1}} \, .
\end{aligned}
\end{equation}
Notice that as the product $a_1 a_n$ grows, the right-hand side of \eqref{varepsilon<0_case} decreases to zero for $\varepsilon > 1$. However, in a similar fashion, $F(\boldsymbol{a})$ grows arbitrarily large as $\min(a_i) \rightarrow \infty$, which yields a contradiction. Thus, it is sufficient to assume that $\varepsilon<1$ holds.

\vspace{2.0mm}

Firstly, let us consider the base case $n=2$. Recall that provided the (weaker) conditions \eqref{conditions on a} hold, then the equality \eqref{Sylvester 2nd Frobenius bound} gives a closed expression for the Frobenius number $F(\boldsymbol{a})$ (likely due to Sylvester \cite{sylvester1884problem}). We prove the following.

\vspace{2.0mm}

\begin{lemma} \label{quadratic_infeasiblity_lemma}
Suppose an integral vector $\boldsymbol{a}=\left(a_1, a_2\right)^T$ satisfies \eqref{conditions on a} with $a_1 \le a_2$. Then, for every $\varepsilon>0$, there is no constant $C > 0$ (independent of $\boldsymbol{a}$) such that 
$$
F(\boldsymbol{a}) \le\ C \cdot \big( a_{1}\,a_{2} \big)^{1-\varepsilon}
$$
holds for all vectors $\boldsymbol{a}$ satisfying \eqref{conditions on a}.
\end{lemma}

\vspace{2.0mm}

\begin{proof}
Suppose for contradiction that there exists a constant $C>0$ and some $\varepsilon>0$ such that for all coprime pairs $(a_1, a_2)$, we have
$$
F(a_1, a_2) \le C \cdot \big( a_{1}\,a_{2} \big)^{1-\varepsilon}
$$

\vspace{2.0mm}

Consider $a_1 = p$ and $a_2 = p+1$, where $p$ denotes some prime. Notice that $\gcd(p, p+1) = 1$. Thus, from the closed expression \eqref{Sylvester 2nd Frobenius bound}, we have
$$
\begin{aligned}
F(p, p+1) &= p \, (p+1) - p - (p+1) \\
            &= p^2 - p - 1.
\end{aligned}
$$

\vspace{2.0mm}

Thus, by assumption, we yield the that inequality
\begin{equation} \label{assumed_contra_base}
F(p, p+1) = p^2 - p - 1 \le  C \cdot \big( p \, (p+1) \big)^{1-\varepsilon}
\end{equation}
holds for every prime $p$. We will show that this inequality \eqref{assumed_contra_base} fails for sufficiently large $p$.

\vspace{2.0mm}

We will use a pair of simple bounds to show the above inequality fails. First, observe the lower bound
$$
p^2 - p - 1 \ge \frac{1}{2} p^2
$$
holds for all $p \ge 2$. Further, note that the upper bound 
$$
p \, (p+1) = p^2 + p \le p^2 + p^2 = 2 p^2,
$$
holds, which yields, upon substitution, the bound
$$
\begin{aligned}
\big( p \, (p+1) \big)^{1-\varepsilon} &\le \big( 2 p^2 \big)^{1-\varepsilon} \\
                                        &= 2^{1-\varepsilon} \, p^{2(1-\varepsilon)}
\end{aligned}
$$
appearing in the right-hand side of \eqref{assumed_contra_base}.

\vspace{2.0mm}

Upon substituting these bounds into \eqref{assumed_contra_base}, we find (for large enough $p$) that
$$
\frac{1}{2} p^2 \le p^2 - p - 1 \le  C \cdot \big( p \, (p+1) \big)^{1-\varepsilon} \le C \cdot \big( 2^{1-\varepsilon} \, p^{2(1-\varepsilon)} \big)
$$
Simplifying further, this yields
\begin{equation} \label{contradiction_yields}
\frac{1}{2} \le  C \cdot 2^{1-\varepsilon} \, p^{2(1-\varepsilon) - 2}.
\end{equation}
However, $p^{2(1-\varepsilon) - 2} = p^{-2\varepsilon}$. Further, observe that as $p \rightarrow \infty$, we have $p^{-2\varepsilon} \rightarrow 0$. Thus, for large $p$, the right-hand side of \eqref{contradiction_yields} becomes arbitrarily small, while the left-hand side remains 1/2. This contradiction completes the proof.
\end{proof}

\vspace{2.0mm}

The results in Tables \ref{tab:test_bounds} and \ref{tab:test_bounds_large_epsilon} provide a concrete demonstration that the inequality \eqref{assumed_contra_base}, which assumes the existence of a universal constant  $C > 0$ fails for sufficiently large primes $p$. Specifically, the tables show that $F(p, p+1) = p^2 - p - 1$  exceeds the bound  $C \cdot \big(p \, (p+1) \big)^{1-\varepsilon}$  for certain choices of $\varepsilon$. It should be noted that in both Table \ref{tab:test_bounds} and \ref{tab:test_bounds_large_epsilon} we fix $C = 1$ for simplicity. The choice of $C$ directly impacts the scale of the test bounds, where larger values of $C$ would proportionally increase the bounds, thereby requiring larger primes $p$ for $F(p, p+1)$ to exceed the test bound. This behaviour aligns with the fact that $F(p, p+1)$ grows quadratically, and the constant $C$ informally {controls} the threshold at which a violation occurs. 

\vspace{2.0mm}

\begin{figure}[ht!]
\centering
	\includegraphics[width=0.8\linewidth]{{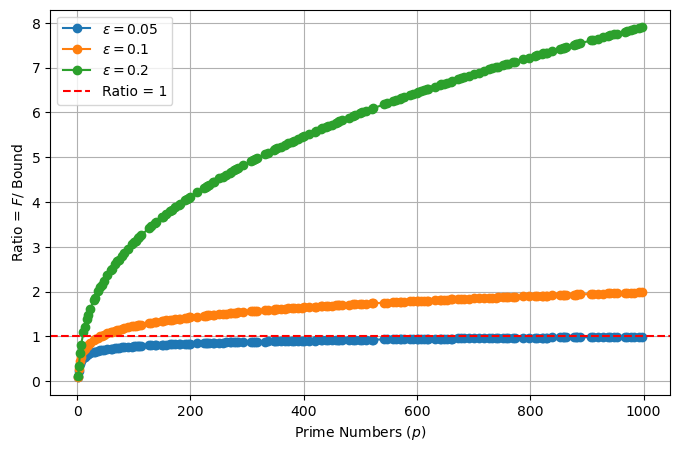}}
	\caption{Ratio of $F(p, p+1) = p^2 - p - 1$ to the test bound $C \cdot \big( p\, (p+1) \big)^{1-\varepsilon}$ for $C=2$ and varying values of $\varepsilon$. Smaller values of $\varepsilon$ produce stricter bounds, keeping the ratio below $1$ for larger primes. Larger $\varepsilon$ relaxes the bound, allowing the ratio to exceed $1$ for smaller primes. The horizontal (red) line at $1$ marks the threshold where $F(p, p+1)$ equals the bound.}
	\label{fig:ratio_vs_primes}
\end{figure}

Figure \ref{fig:ratio_vs_primes} illustrates the ratio 
$$
\frac{F(p, p+1)} { C \cdot \big(p \, (p+1)\big)^{1-\varepsilon}}
$$ 
plotted against prime values $p$ for different values of $\varepsilon$. Each curve corresponds to a different value of $\varepsilon$, but using the same constant $C=2$. The horizontal (red) line at $1$ indicates the exact point where $F(p,p+1)$ and the bounding function coincide. Points above this line reflect primes $p$ for which 
$F(p,p+1) > C\,\big(p\,(p+1)\big)^{1-\varepsilon}$ 
holds, thus violating the proposed bound. 
Observe that smaller values of $\varepsilon$ lead to a \enquote{stricter} (larger) exponent $(1-\varepsilon)$ and thus a larger sub-quadratic bound, causing the ratio to remain below $1$ for more primes. Conversely, {larger} $\varepsilon$ yields a more \enquote{relaxed} (smaller) exponent, so $F(p,p+1)$ exceeds the bound {earlier} (i.e. for smaller primes). Overall, this illustrates that for each fixed $\varepsilon>0$, there exists a large enough prime $p$ such that $F(p,p+1)$ ultimately outgrows $C \cdot \big( p \, (p+1)\big)^{\,1-\varepsilon}$.

\vspace{2.0mm}

It should be emphasised that in Lemma \ref{quadratic_infeasiblity_lemma} we have established that no sub‐$\big(a_1 \cdot a_2\big)$ bound can hold in dimension $n=2$.  In particular, choosing $a_1=p$ and $a_2=p+1$ shows that
$$
F(p,p+1) \;=\; p^2 - p - 1
\quad
\text{eventually outgrows }
C \cdot \bigl(p\,(p+1)\bigr)^{1-\varepsilon}
\text{ for any fixed }C>0 \text{ and } \varepsilon>0.
$$
Thus, the theorem statement is already proven in the two-dimensional setting.  Next, we handle the case that $n\ge3$ using a standard {embedding} argument.

\vspace{2.0mm}

\begin{proof}
Suppose for contradiction that there exists some integer $n \ge 3$, a constant $C>0$ and some $\varepsilon>0$ such that 
$$
F(\boldsymbol{a}) \le C \cdot \bigl(a_1\,a_n\bigr)^{1-\varepsilon}
$$
holds for every integer vector $\boldsymbol{a} \in \mathbb{Z}^n_{>0}$ satisfying \eqref{conditions on a}. 

\vspace{2.0mm}

We extend the \enquote{problematic} two-dimensional pair $(p,p+1)$ from Lemma \ref{quadratic_infeasiblity_lemma} into an $n$-dimensional vector for $n \ge 3$, defined as
$$
\widetilde{\boldsymbol{a}}
=
\Bigl(p,\;p+1,\;\;\underbrace{r,\,r,\dots,r}_{n-2\text{ times}}\Bigr)^T,
$$
where $r \in \mathbb{Z}_{>0}$ such that $r \ge p+1$ is \enquote{carefully} chosen. In particular, we aim to ensure that $\gcd\left(\widetilde{\boldsymbol{a}}\right) = 1$, that $r$ is representable by the pair $\{p,p+1\}$ (i.e. $r = px + (p+1)y$ holds for some $x,y \in \mathbb{Z}_{\ge 0}$) and that $r \approx cp$ for some constant $c>1$. 

\vspace{2.0mm}

Observe that $\gcd\left(\widetilde{\boldsymbol{a}}\right) = 1$ is guaranteed since $\gcd(p, p+1) = 1$ for all primes $p$. Let us pick some $(x_0, y_0)$ such that 
$$
r = p \, x_0 + (p+1) \, y_0 \ge p+1,
$$
which ensures that the entries of $\widetilde{\boldsymbol{a}}$ are nondecreasing and that $r$ is representable by the pair $\{p,p+1\}$ and that $r \approx cp$ for some constant $c>1$. 

\vspace{2.0mm}

Moreover, by construction, we have
$$
F(\widetilde{\boldsymbol{a}}) = F\big(p,\;\;p+1,\;\;\underbrace{r,\,r,\dots,r}_{n-2\text{ times}}\big) = F(p, p+1) = p^2 - p - 1,
$$
which follows since clearly no integer $m \le p^2 - p - 1$ which is unrepresentable by $\{p,p+1\}$ becomes representable by $\{p,p+1,r\}$.

\vspace{2.0mm}

Thus, by assumption, we deduce that
\begin{equation} \label{assumed_contra_nge3}
\begin{aligned}
F(\widetilde{\boldsymbol{a}}) = p^2 - p - 1 &\le C \cdot \bigl(a_1\,a_n\bigr)^{1-\varepsilon} \\
                                &= C \cdot \bigl(p\,r\bigr)^{1-\varepsilon}
\end{aligned}
\end{equation}
holds for every prime $p$. We will show that this inequality \eqref{assumed_contra_nge3} fails for sufficiently large $p$.

\vspace{2.0mm}

Recall that $r \approx cp$, hence \eqref{assumed_contra_nge3} yields
\begin{equation} \label{assumed_contra_nge3111}
\begin{aligned}
p^2 - p - 1 &\le C \cdot \bigl(p \, r\bigr)^{1-\varepsilon} \\
            &\approx C \cdot \bigl( c p^2 \, \bigr)^{1-\varepsilon} \\
            &= C \cdot \big( c^{1-\varepsilon} p^{2(1-\varepsilon)} \big)
\end{aligned}
\end{equation}

\vspace{2.0mm}

Recall that the lower bound 
$$
p^2 - p - 1 \ge \frac{1}{2} p^2
$$ 
holds for all $p \ge 2$. Thus, \eqref{assumed_contra_nge3111} becomes
$$
\frac{1}{2} p^2 \le C \cdot \big( c^{1-\varepsilon} p^{2(1-\varepsilon)} \big) ,
$$
where simplifying further implies
\begin{equation} \label{assumed_contra_nge3111final}
\frac{1}{2} \le C \cdot c^{1-\varepsilon} p^{2(1-\varepsilon) - 2}.
\end{equation}
However, $p^{2(1-\varepsilon) - 2} = p^{-2\varepsilon}$. Further, similarly observe that as $p \rightarrow \infty$, we have $p^{-2\varepsilon} \rightarrow 0$. Thus, for large primes $p$, the right-hand side of \eqref{assumed_contra_nge3111final} becomes arbitrarily small, while the left-hand side remains 1/2. This contradiction completes the proof.
\end{proof}

\vspace{2.0mm}


\begin{table}[ht!]
\centering
\resizebox{\textwidth}{!}{
\begin{tabular}{c | c | c c | c c | c c}
$p$ & $F(p, p+1)$ & Bound ($\varepsilon=0.005$) & $F > \text{Bound}$ & Bound ($\varepsilon=0.01$) & $F > \text{Bound}$ & Bound ($\varepsilon=0.02$) & $F > \text{Bound}$ \\
\hline
    2 &     1 &          5.95 &     False &          5.89 &     False &          5.79 &     False \\
    3 &     5 &         11.85 &     False &         11.71 &     False &         11.42 &     False \\
    5 &    19 &         29.49 &     False &         29.00 &     False &         28.03 &     False \\
    7 &    41 &         54.88 &     False &         53.79 &     False &         51.67 &     False \\
   11 &   109 &        128.82 &     False &        125.71 &     False &        119.72 &     False \\
   13 &   155 &        177.33 &     False &        172.77 &     False &        164.01 &     False \\
   17 &   271 &        297.37 &     False &        288.98 &     False &        272.90 &     False \\
   19 &   341 &        368.88 &     False &        358.08 &     False &        337.43 &      True \\
   23 &   505 &        534.85 &     False &        518.23 &     False &        486.52 &      True \\
   29 &   811 &        841.05 &     False &        813.06 &     False &        759.85 &      True \\
   31 &   929 &        958.36 &     False &        925.86 &      True &        864.13 &      True \\
   37 & 1,331 &      1,355.96 &     False &      1,307.69 &      True &      1,216.26 &      True \\
   41 & 1,639 &      1,659.03 &     False &      1,598.35 &      True &      1,483.59 &      True \\
   43 & 1,805 &      1,821.95 &     False &      1,754.49 &      True &      1,626.98 &      True \\
   47 & 2,161 &      2,170.56 &     False &      2,088.36 &      True &      1,933.18 &      True \\
   53 & 2,755 &      2,750.34 &      True &      2,643.04 &      True &      2,440.82 &      True \\
   59 & 3,421 &      3,398.27 &      True &      3,262.22 &      True &      3,006.24 &      True \\
\end{tabular}}
\caption{Comparison of $F(p, p+1) = p^2 - p - 1$ with test bounds of the form $C \cdot (p(p+1))^{1-\varepsilon}$ for various (smaller) values of $\varepsilon$ with $C=1$. It should be noted that smaller values of $\varepsilon$ test stricter sub-quadratic growth.}
\label{tab:test_bounds}
\end{table}

\begin{table}[ht!]
\centering
\resizebox{\textwidth}{!}{
\begin{tabular}{c | c | c c | c c | c c}
$p$ & $F(p, p+1)$ & Bound ($\varepsilon=0.05$) & $F > \text{Bound}$ & Bound ($\varepsilon=0.1$) & $F > \text{Bound}$ & Bound ($\varepsilon=0.2$) & $F > \text{Bound}$ \\
\hline
    2 &        1 &          5.49 &     False &          5.02 &     False &          4.19 &     False \\
    3 &        5 &         10.60 &     False &          9.36 &     False &          7.30 &     False \\
    5 &       19 &         25.31 &     False &         21.35 &     False &         15.19 &      True \\
    7 &       41 &         45.79 &     False &         37.44 &      True &         25.04 &      True \\
   11 &      109 &        103.41 &      True &         81.01 &      True &         49.71 &      True \\
   13 &      155 &        140.30 &      True &        108.16 &      True &         64.28 &      True \\
\end{tabular}}
\caption{Comparison of $F(p, p+1) = p^2 - p - 1$ with test bounds of the form $C \cdot (p\,(p+1))^{1-\varepsilon}$ for various (larger) values of $\varepsilon$ with $C=1$. These values test less strict sub-quadratic growth.}
\label{tab:test_bounds_large_epsilon}
\end{table}

\section{Conclusions and Future Work}
In this paper, we conducted a detailed comparative study of known upper bounds on the Frobenius number $F(\boldsymbol{a})$, focusing on both classical and recently proposed bound. Through Monte Carlo simulations, we evaluated the empirical tightness of each bound under two distinct assumptions on the input $\boldsymbol{a}$, namely the (weaker) conditions \eqref{conditions on a} and the (stronger) pairwise coprimality conditions \eqref{stronger conditions on a}.

\vspace{2.0mm}

Under the (weaker) conditions \eqref{conditions on a}, classical bounds such as those of Erd{\H o}s and Graham, Schur, and Vitek typically perform well, with moderate errors and predictable behaviour. The bound proposed by Fukshansky and Robins, in contrast, exhibits significantly higher errors, particularly as dimension and the maximum absolute entry of $\boldsymbol{a}$ increase, limiting its practical applicability. Under the (stronger) pairwise coprimality conditions \eqref{stronger conditions on a}, the Sylvester-based extension due to Williams and Haijima performs well, although suffers from higher variability in extreme cases. The bound of Selmer performs well when the maximum absolute entry of $\boldsymbol{a}$ is large and the minimal absolute entry of $\boldsymbol{a}$ is moderate. 

\vspace{2.0mm}

We additionally established a key theoretical result (Theorem \ref{quadratic_infeasiblity}), showing that under the (weaker) conditions \eqref{conditions on a}, no general upper bound on the Frobenius number can exhibit sub-quadratic growth. This builds upon the recent work of Williams and Haijima \cite{williams2023considering} and confirms that quadratic dependence on the largest entry is unavoidable unless additional assumptions are imposed.

\vspace{2.0mm}

There are several promising directions for future research. Firstly, developing hybrid or adaptive bounding techniques, informed by empirical findings, would allow dynamic selection of the most suitable bound based on structural properties of the input, potentially offering both reliability and sharpness. Secondly, the limitations proven in this paper motivate further theoretical investigations into conditions under which sub-quadratic bounds might be possible. For instance, studying restricted classes of input vectors characterised by sparsity, bounded entry ratios $a_n/a_1$, or lattice-width constraints could reveal scenarios in which tighter bounds are feasible. Thirdly, future studies could explore probabilistic bounding approaches, shifting from deterministic worst-case guarantees to high-probability bounds for random input vector distributions. Such probabilistic analyses could potentially yield tighter practical bounds for typical-case scenarios, even when deterministic sub-quadratic bounds remain unattainable. Fourthly, embedding these theoretical bounds as heuristics or pruning criteria within exact Frobenius number algorithms is an important computational direction. Fifthly, data-driven and machine learning methods represent another promising avenue. Future research could employ regression or kernel-based learning models trained on simulated data to predict $F(\boldsymbol{a})$ directly. Insights from these models could subsequently inform the theoretical refinement of analytic bounds, revealing new structural patterns that classical analytical methods alone may not detect. Finally, expanding this comparative framework to modular Frobenius numbers, weighted Frobenius problems, or coin sets forming arithmetic progressions, would help extend the relevance and applicability of these findings to a wider class of integer optimisation problems.

\newpage





\newpage
\appendix
\section{Additional Box Plots for Larger Inputs} \label{appendix:larger_m}
The following demonstrates that with increases in the upper bound $m$, the scale of the errors grows, particularly for bounds involving the $\ell_2$-norm and gamma function. These figures further support the conclusion that certain bounds become less practical at large input sizes.

\vspace{2.0mm}

\begin{figure}[ht!]
    \centering
    \includegraphics[width=1\linewidth]{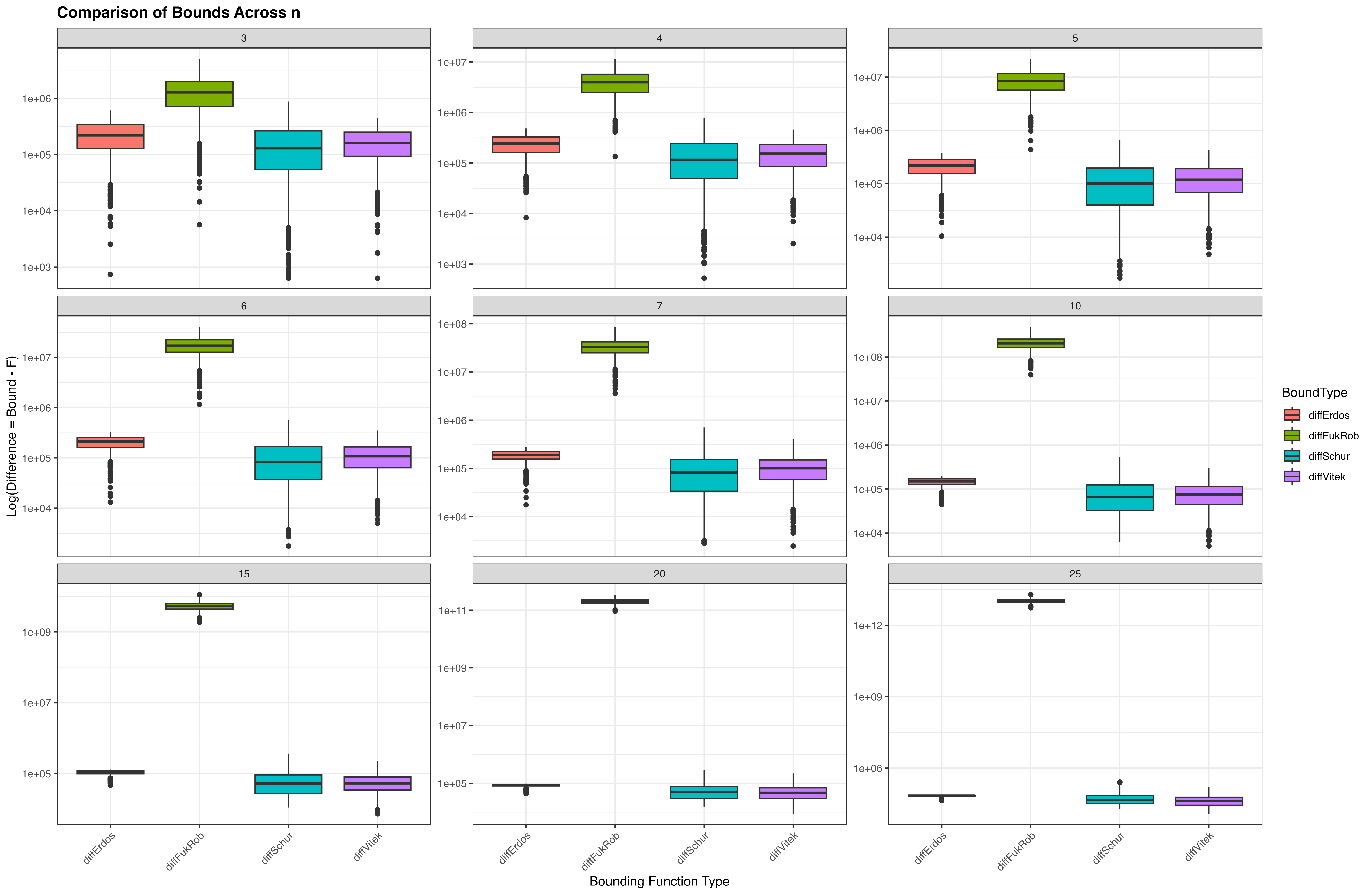}
    \caption{Box plots of the difference, i.e. \(\text{Bound} - F(\boldsymbol{a})\), with \(\|\boldsymbol{a}\|_\infty \le 1000\).}
    \label{fig:boxplot_appendix_m1000}
\end{figure}

\newpage

\begin{figure}[ht!]
    \centering
    \includegraphics[width=1\linewidth]{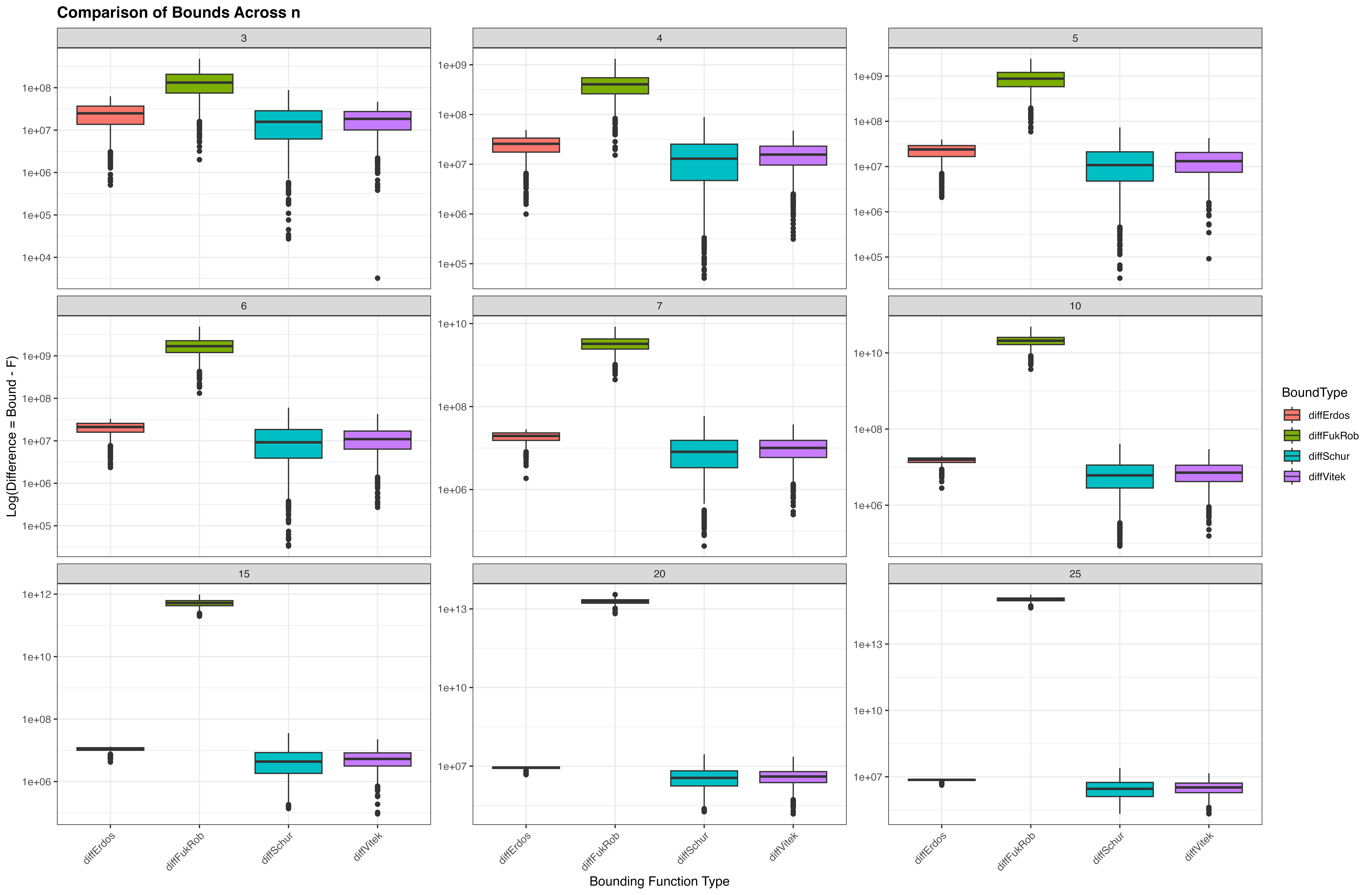}
    \caption{Box plots of the difference, i.e. \(\text{Bound} - F(\boldsymbol{a})\), with \(\|\boldsymbol{a}\|_\infty \le 10{,}000\).}
    \label{fig:boxplot_appendix_m10000}
\end{figure}

\newpage
\section{Box plots with Sample Density for Larger Inputs} \label{appendix:densities}
The following figures extend the analysis in the main text to much larger input scales. Note that as teh upper bound \( m \) increases, the overestimation of certain bounds (notably Fukshansky and Robins) becomes more dramatic, consistent with the bound’s theoretical dependence on norms and gamma functions.

\vspace{2.0mm}

\begin{figure}[ht!]
    \centering
    \includegraphics[width=0.98\linewidth]{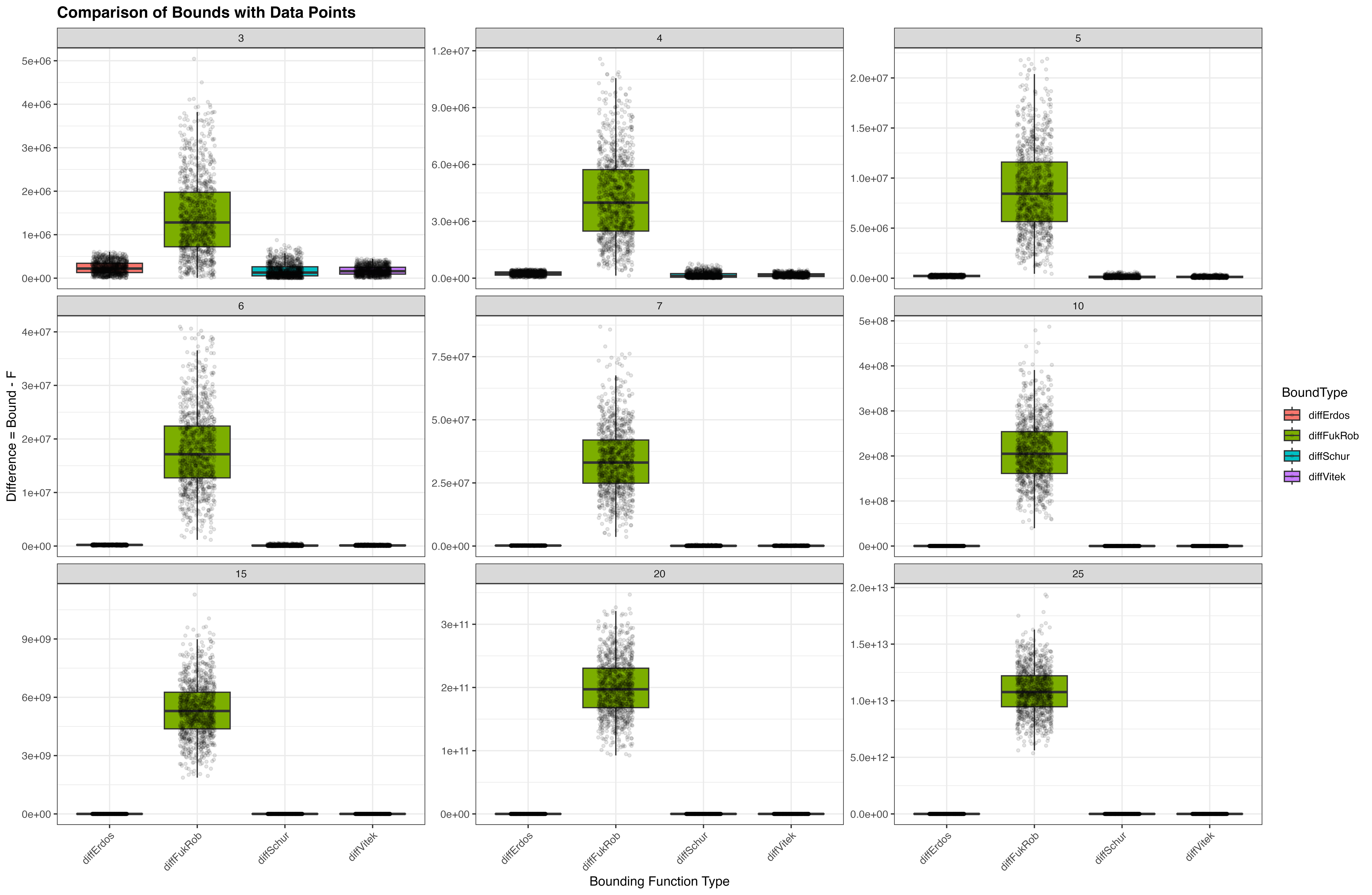}
    \caption{Box plots with overlaid sample points for \(\|\boldsymbol{a}\|_\infty \le 1000\).}
    \label{fig:boxplot_data_m1000}
\end{figure}

\newpage
\begin{figure}[ht!]
    \centering
    \includegraphics[width=1\linewidth]{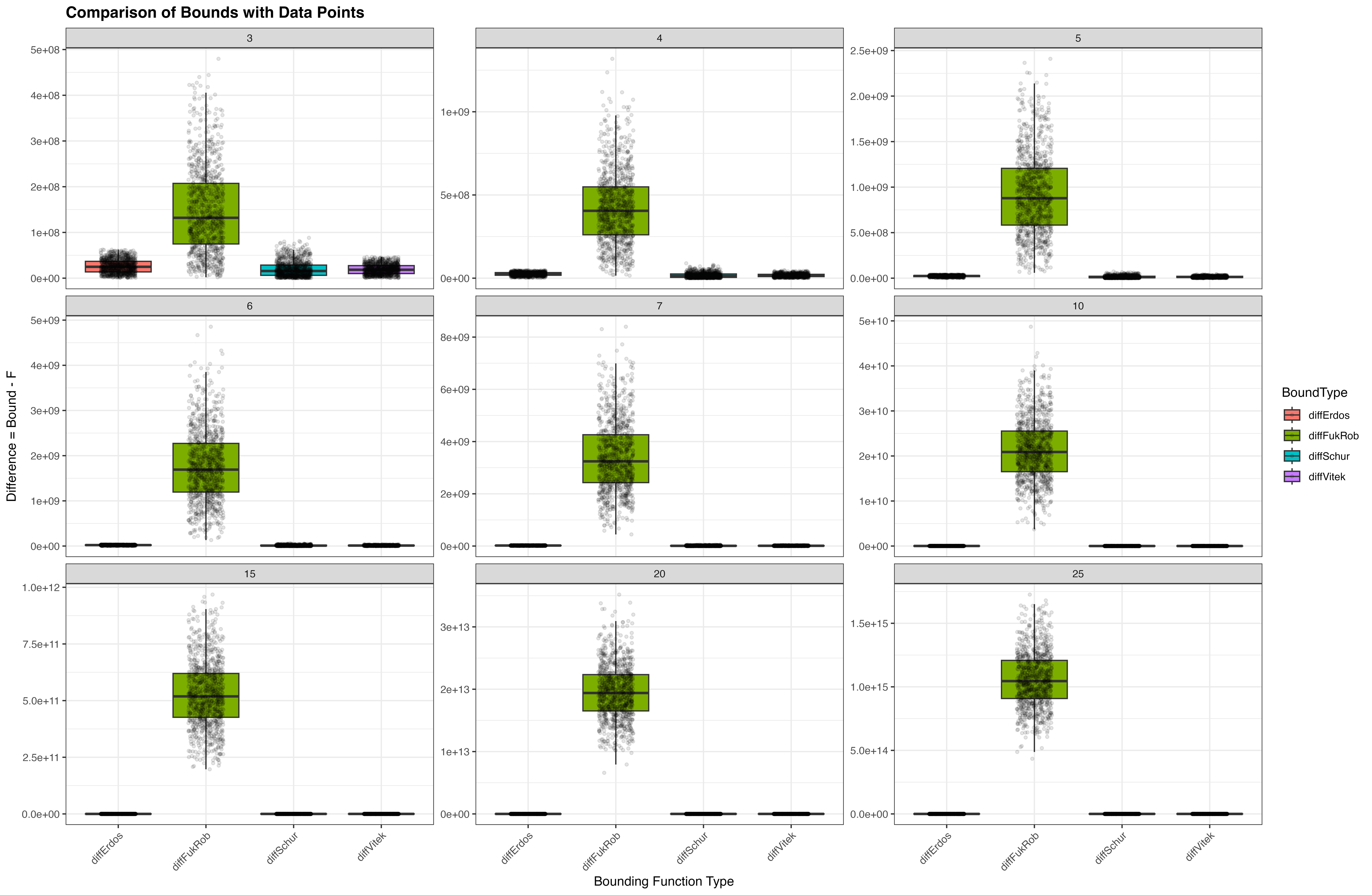}
    \caption{Box plots with overlaid sample points for \(\|\boldsymbol{a}\|_\infty \le 10{,}000\).}
    \label{fig:boxplot_data_m10000}
\end{figure}

\newpage
\section{Density Plots for Larger Inputs}
\label{appendix:density_plots}
The following density plots reinforce the findings from the main text. Observe that as the bound on the maximum entry increases, the distribution of errors for most classical bounds remain concentrated, while the Fukshansky and Robins bound displays increasingly heavy tails and long-range deviations. This empirical behaviour aligns with its theoretical scaling structure.

\vspace{2.0mm}

\begin{figure}[ht!]
    \centering
    \includegraphics[width=0.9\linewidth]{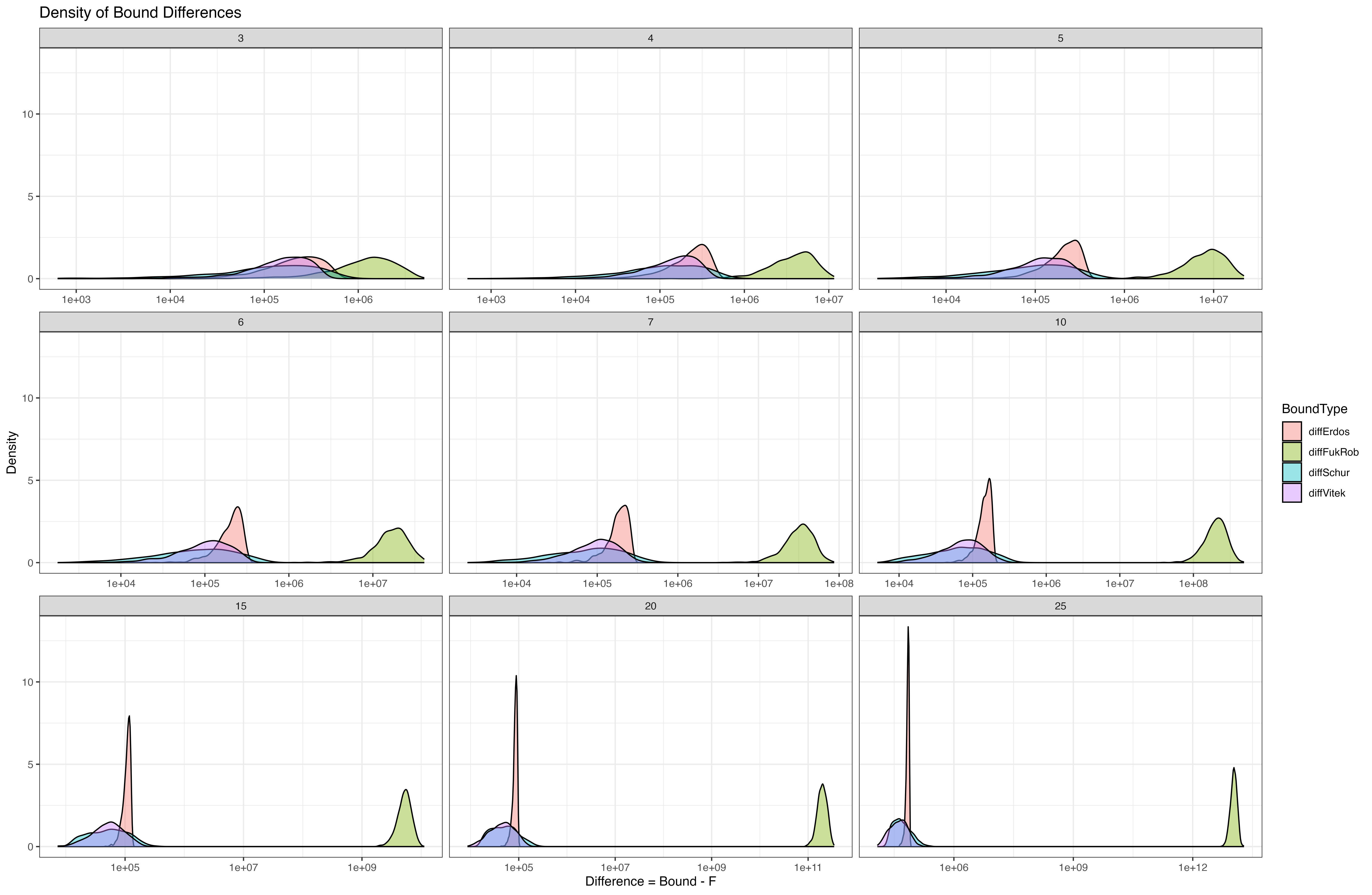}
    \caption{Density plots of the difference, i.e. \(\text{Bound} - F(\boldsymbol{a})\), for \(\|\boldsymbol{a}\|_\infty \le 1000\).}
    \label{fig:density_m1000}
\end{figure}

\newpage
\begin{figure}[ht!]
    \centering
    \includegraphics[width=1\linewidth]{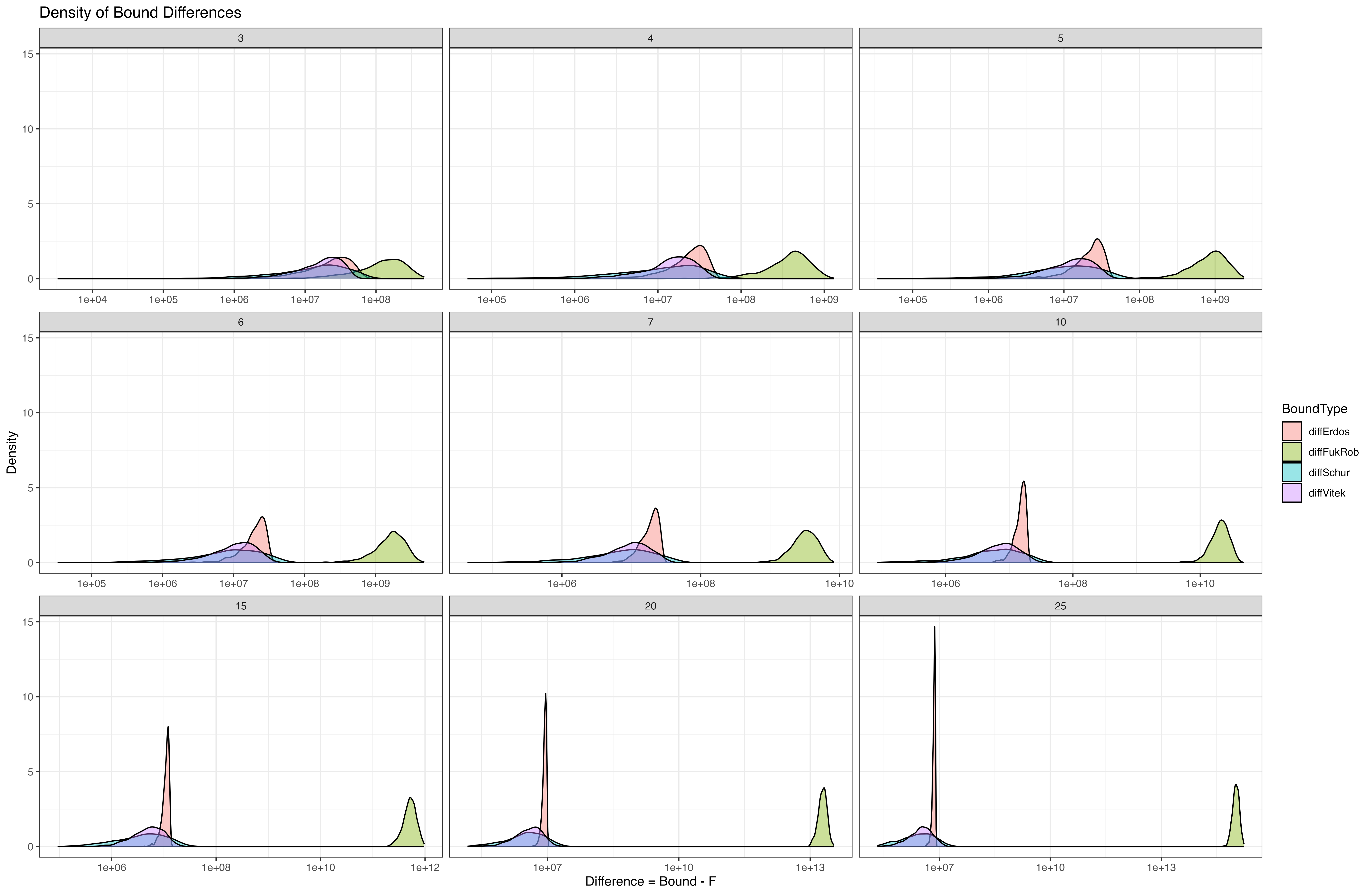}
    \caption{Density plots of the difference, i.e. \(\text{Bound} - F(\boldsymbol{a})\), for \(\|\boldsymbol{a}\|_\infty \le 10{,}000\).}
    \label{fig:density_m10000}
\end{figure}

\newpage
\section{Line Plots of Mean Error vs. $n$ for Larger Inputs} \label{appendix:line_error_plots}
The following line plots for larger values of the upper bound $m$ illustrate the diverging behaviour between classical and some modern bounds. Observe that while the bounds seem to increase with the dimension $n$, the rate at which the Fukshansky and Robins bound in expectation becomes weaker is markedly steeper, confirming its theoretical super-exponential dependence.

\vspace{2.0mm}

\begin{figure}[ht!]
    \centering
    \includegraphics[width=0.95\linewidth]{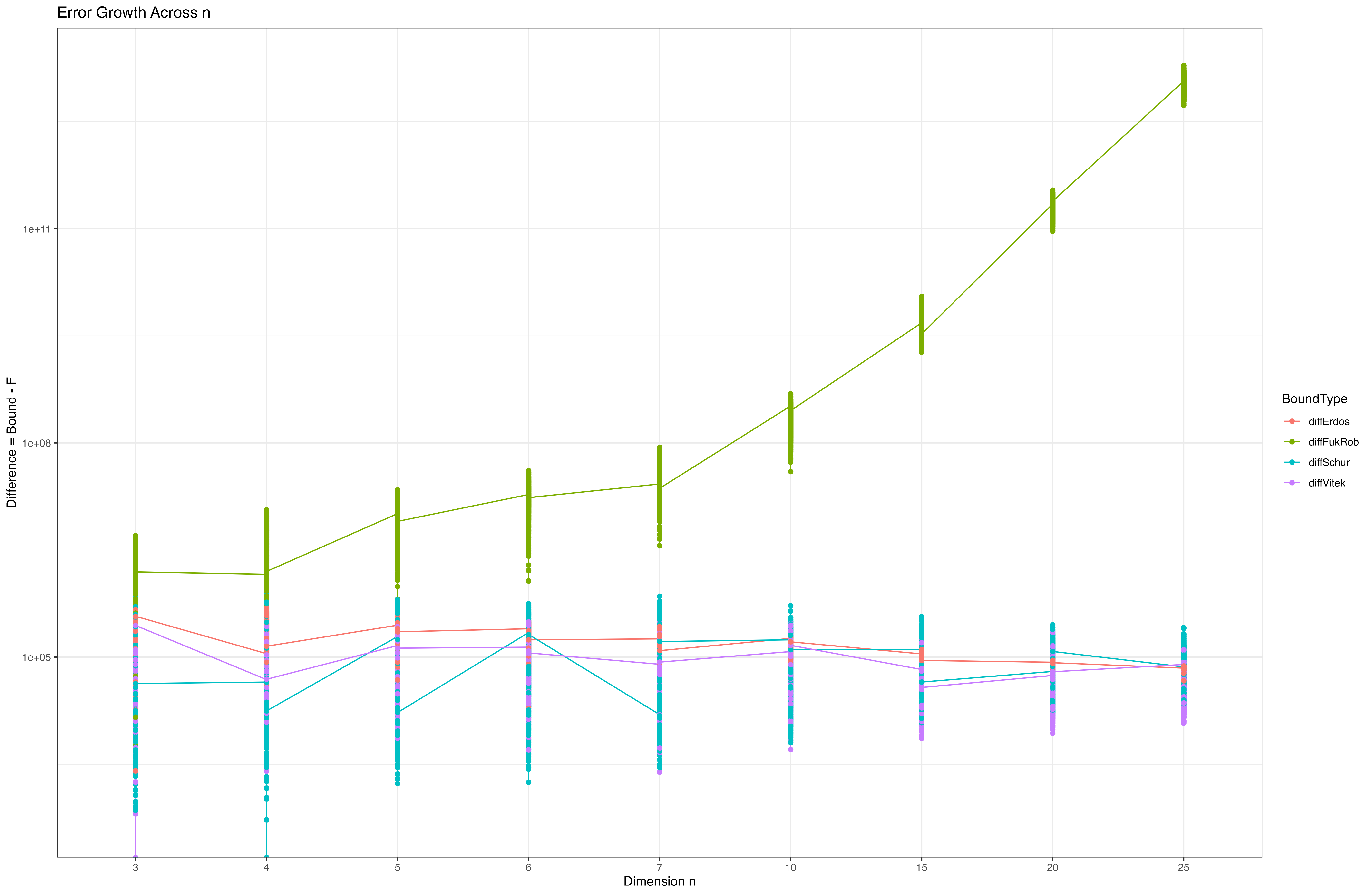}
    \caption{Line plots of the difference, i.e. \(\text{Bound} - F(\boldsymbol{a})\), plotted against the dimension \(n\) for vectors with \(\|\boldsymbol{a}\|_\infty \le 1000\).}
    \label{fig:line_error_m1000}
\end{figure}

\newpage
\begin{figure}[ht!]
    \centering
    \includegraphics[width=1\linewidth]{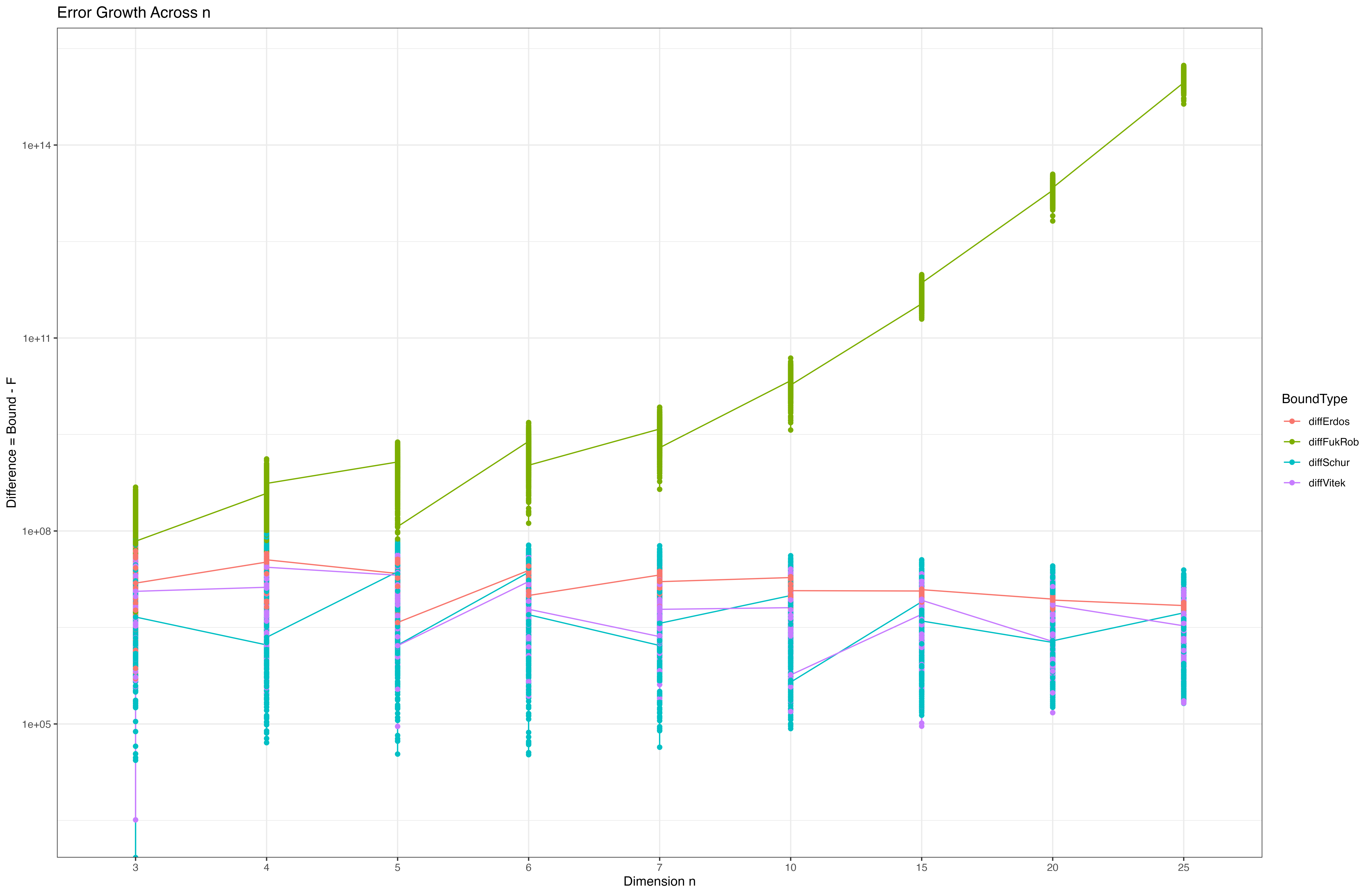}
    \caption{Linear plots of the difference, i.e. \(\text{Bound} - F(\boldsymbol{a})\), plotted against the dimension \(n\) for vectors with \(\|\boldsymbol{a}\|_\infty \le 10{,}000\). }
    \label{fig:line_error_m10000}
\end{figure}

\newpage
\section{Scatter Plots of Fukshansky and Robins Error vs. Maximum Entry}
\label{appendix:scatter_plots}
The following scatter plots provide additional evidence of the scaling behaviour of the Fukshansky and Robins bound. It should be emphasised that while other bounds remain somewhat insensitive to growth in the maximum input entry, the error associated with the bound of Fukshansky and Robins increases steeply as $\|\boldsymbol{a}\|_\infty$ grows. This pattern is exacerbated by dimension, affirming the theoretical structure involving high-order norm and gamma dependencies.

\vspace{2.0mm}

\begin{figure}[ht!]
    \centering
    \includegraphics[width=1\linewidth]{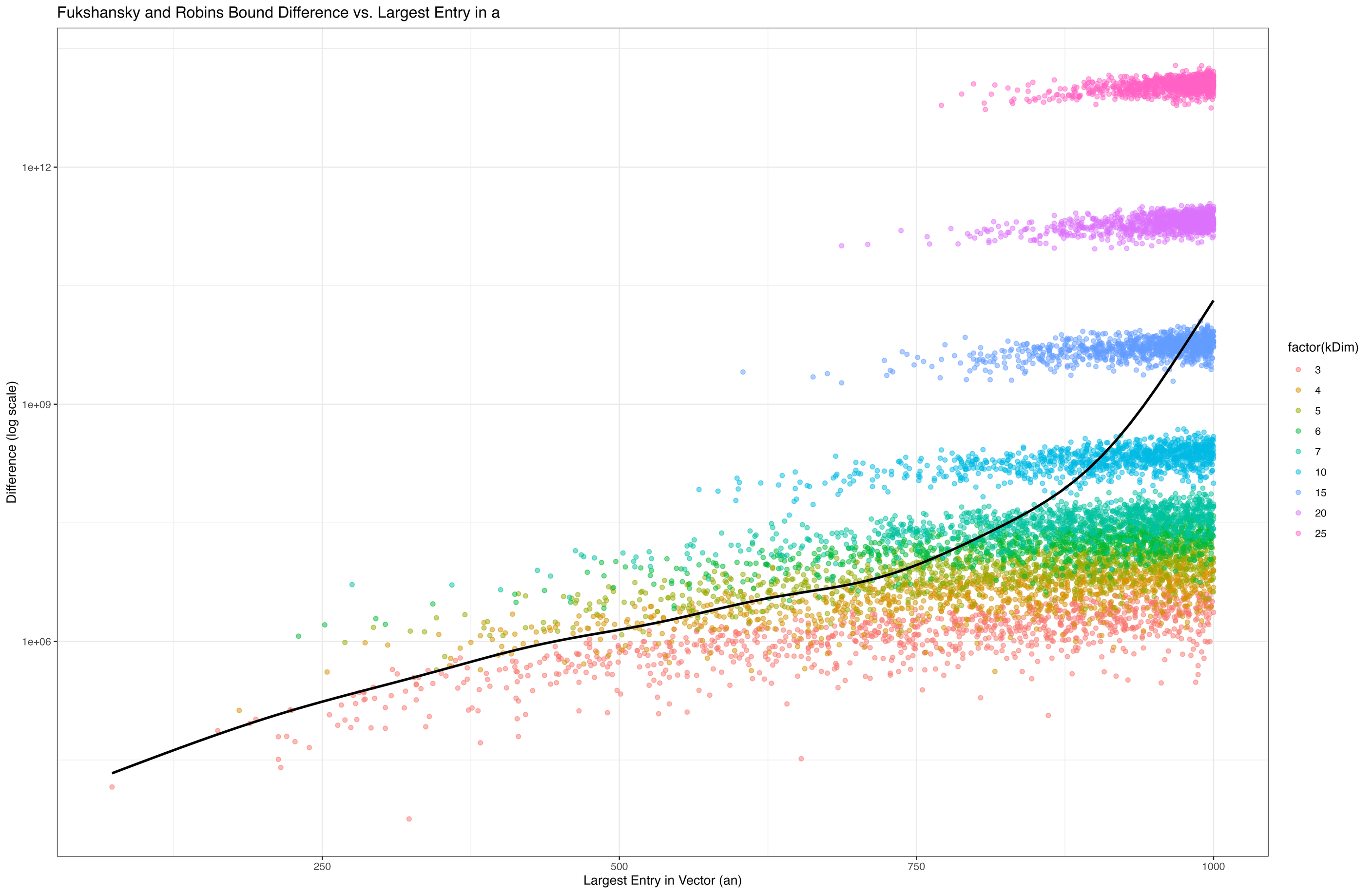}
    \caption{Scatter plot of the difference \(\texttt{FukRob} - F(\boldsymbol{a})\) plotted against \(\max_i a_i\) for vectors with \(\|\boldsymbol{a}\|_\infty \le 1000\). }
    \label{fig:scatter_fukrob_m1000}
\end{figure}

\newpage
\begin{figure}[ht!]
    \centering
    \includegraphics[width=1\linewidth]{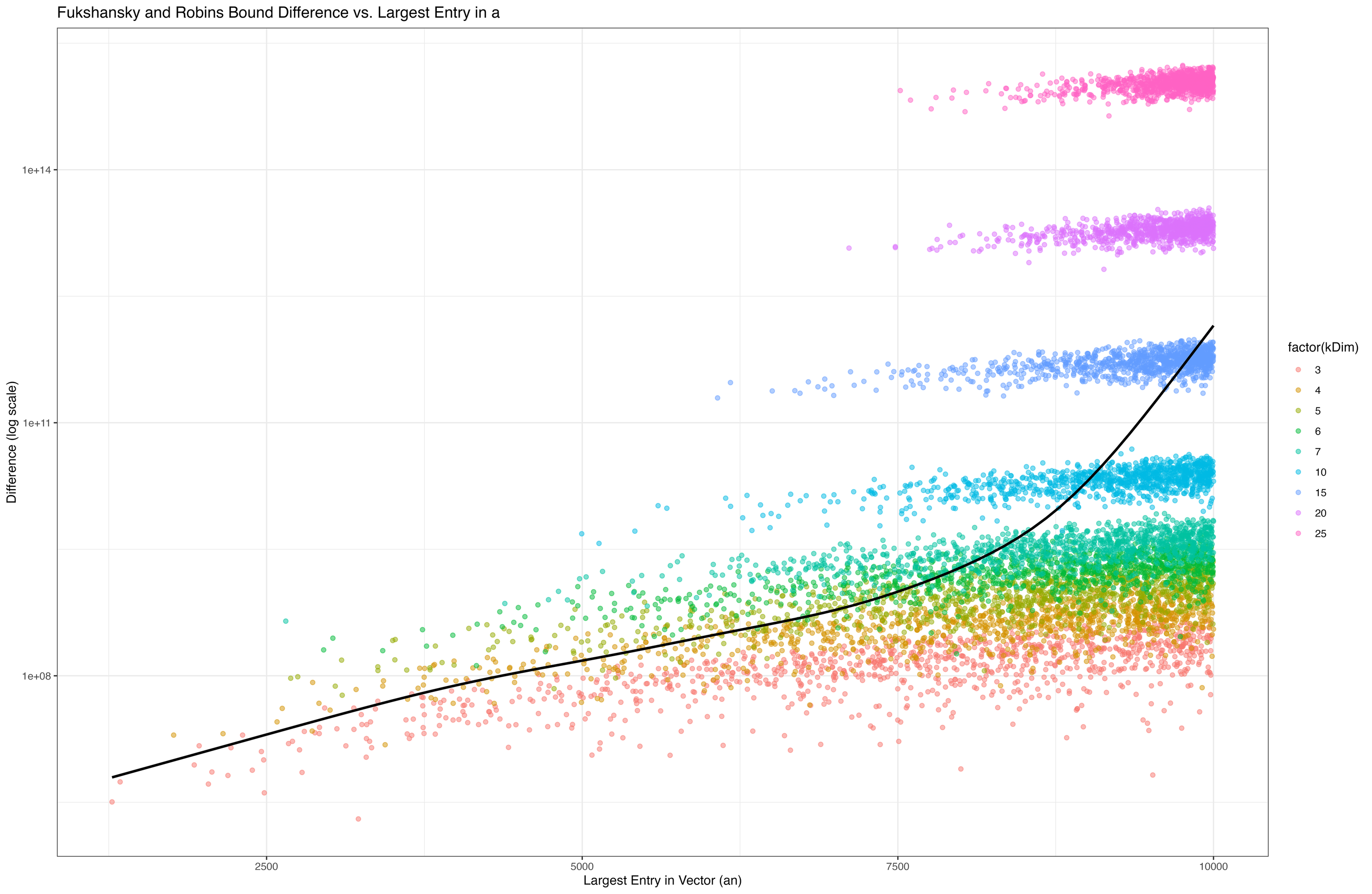}
    \caption{Scatter plot of \(\texttt{FukRob} - F(\boldsymbol{a})\) plotted against \(\max_i a_i\) for vectors with \(\|\boldsymbol{a}\|_\infty \le 10{,}000\). }
    \label{fig:scatter_fukrob_m10000}
\end{figure}

\newpage
\section{Box Plots of Bound Errors} \label{appendix:box_m10000}



The following figure complements Figure~\ref{fig:boxplot_all_k}, with inputs restricted to $\|\boldsymbol{a}\|_\infty \le 10{,}000$. The overall distributional shapes are consistent, but we observe reduced vertical spread and fewer extreme outliers across all bounds. This supports the view that large-magnitude entries drive some erratic behaviour in bound performance.

\vspace{2.0mm}

\begin{figure}[ht!]
\centering
\includegraphics[width=1\linewidth]{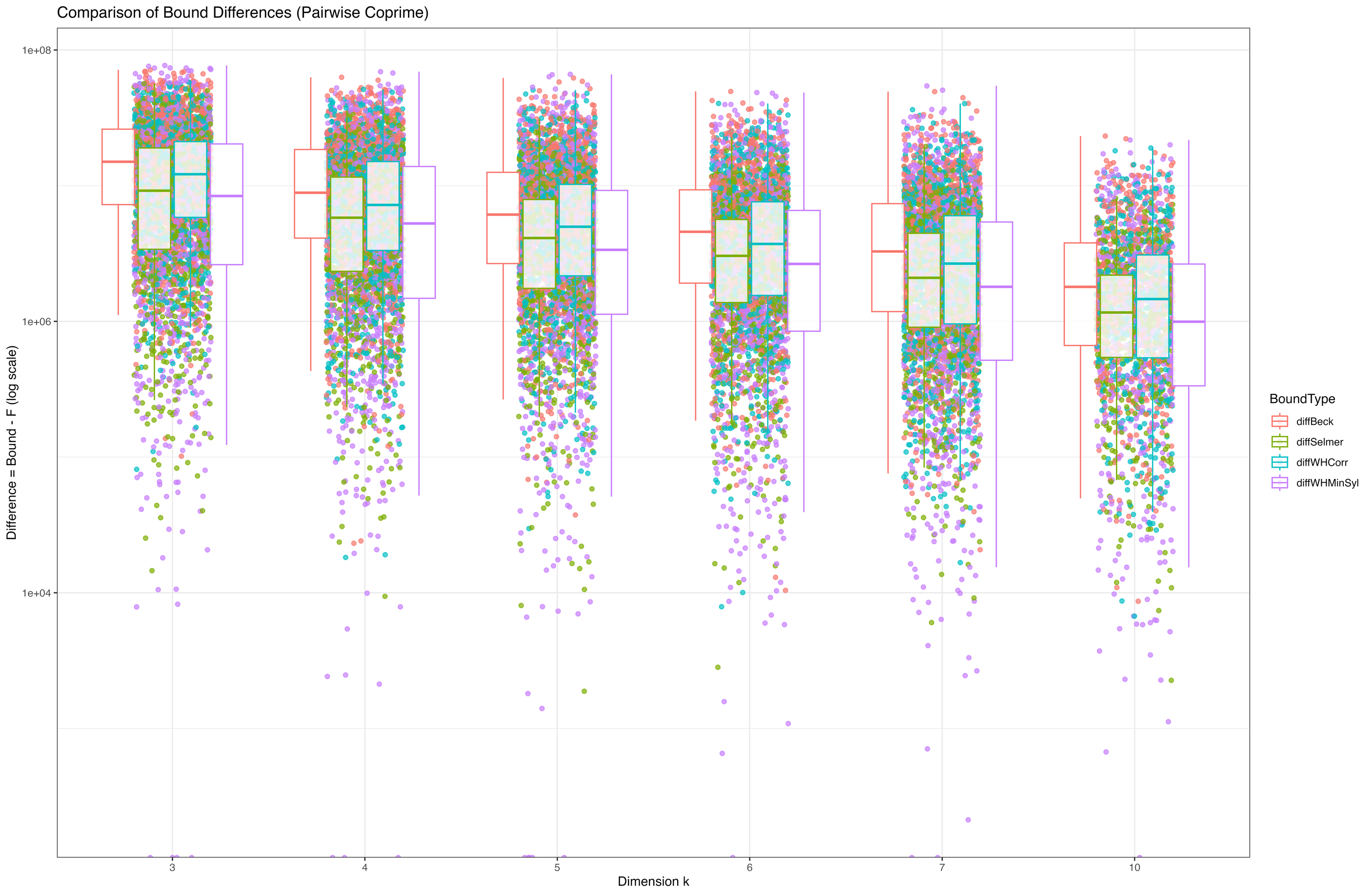}
\caption{Comparison of the difference, i.e. \(\text{Bound} - F(\boldsymbol{a})\), across bounds with \(\|\boldsymbol{a}\|_\infty \le 10{,}000\).}
\label{fig:scatter_largest_an_m10000}
\end{figure}

\newpage
\section{Kernel Density Estimates Grouped by Dimension and Bound Type} 
\label{appendix:density_structured_m10000}

The following figure complements Figure \ref{fig:density_by_k} in the main text, but now with $\| \boldsymbol{a} \|_{\infty} \le 10,000$. 

\vspace{2.0mm}

\begin{figure}[ht!]
\centering
\includegraphics[width=1\linewidth]{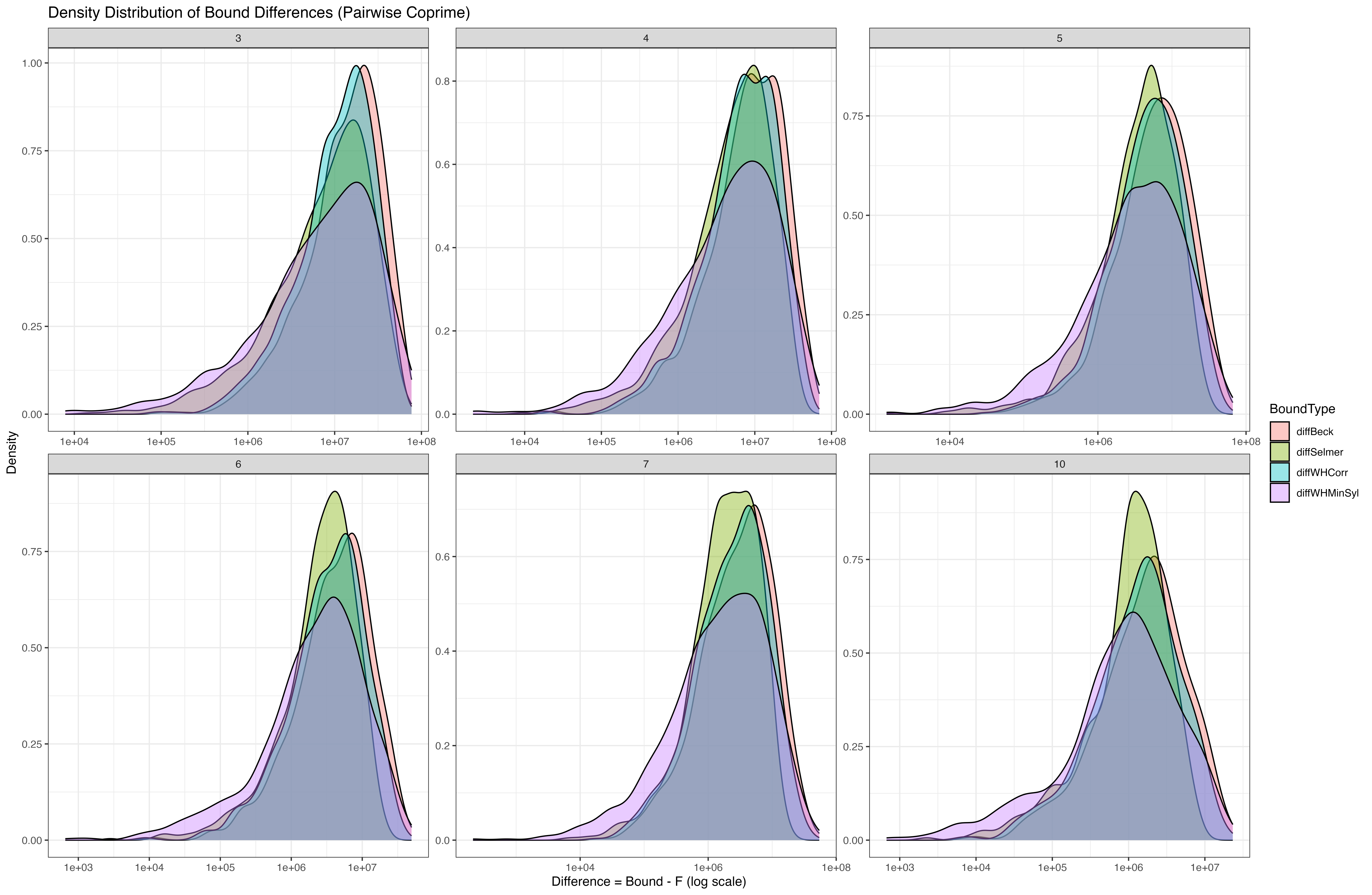}
\caption{Density plots of the difference, i.e. $\text{Bound} - F(\boldsymbol{a})$, for $\| \boldsymbol{a} \|_{\infty} \le 10,000$ across dimensions $n$. }
\label{fig:density_m10000_structured}
\end{figure}

\newpage
\section{Scatter Plot of Bound Differences Grouped by Type} 
\label{appendix:boxplot_m10000_structured}

The following figure compliments Figure~\ref{fig:scatter_diff_vs_an} in the main text, but with inputs satisfying $\| \boldsymbol{a} \| \le 10,000$. The black curve is a smoothed generalized additive model (GAM) curve, which highlights the general scaling trend.

\vspace{2.0mm}

\begin{figure}[ht!]
\centering
\includegraphics[width=1\linewidth]{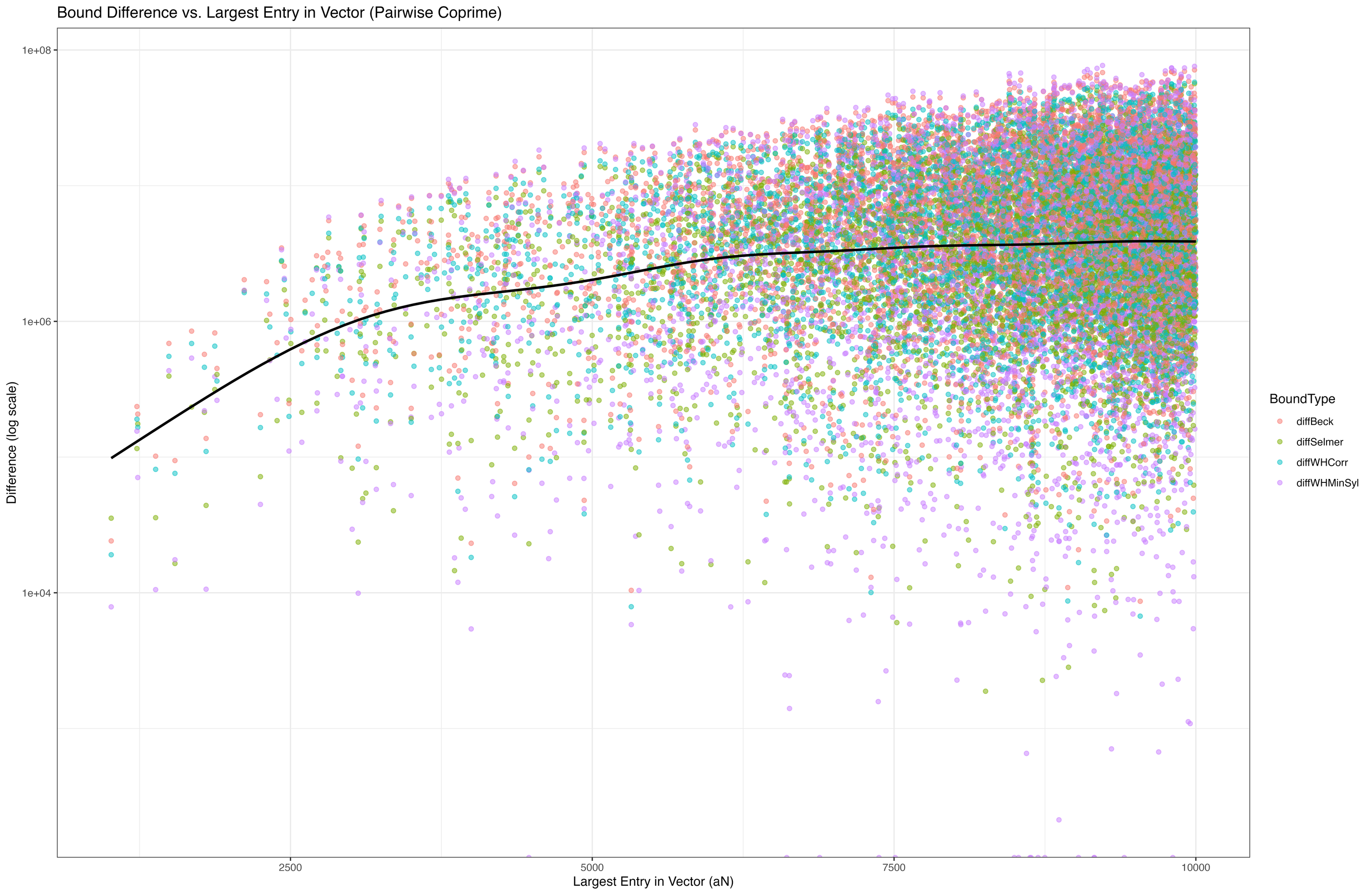}
\caption{Scatter plot of the difference, i.e. \(\text{Bound} - F(\boldsymbol{a})\), plotted against the maximum entry $a_n$, coloured by bound type, with GAM smoother for $\| \boldsymbol{a} \|_{\infty} \le 10,000$.}
\label{fig:boxplot_structured_m10000}
\end{figure}

\newpage
\section{Relative Performance Regions}
\label{appendix:best_bound_region_m10000}

The following plot complements Figure~\ref{fig:scatter_best_bound_by_inputs} in the main text. Note that each point corresponds to an input vector $\boldsymbol{a}$, which is then coloured according to the bound that achieved the lowest error for that random instance. Observe that the second bound of Williams and Haijima \eqref{williams_sylvester_extension}
performs particularly well in regions with small $a_1$, while the upper bound of Selmer \eqref{Selmer_bound_n_1} often achieves the lowest error when $a_1$ grows.

\vspace{2.0mm}

\begin{figure}[ht!]
\centering
\includegraphics[width=1\linewidth]{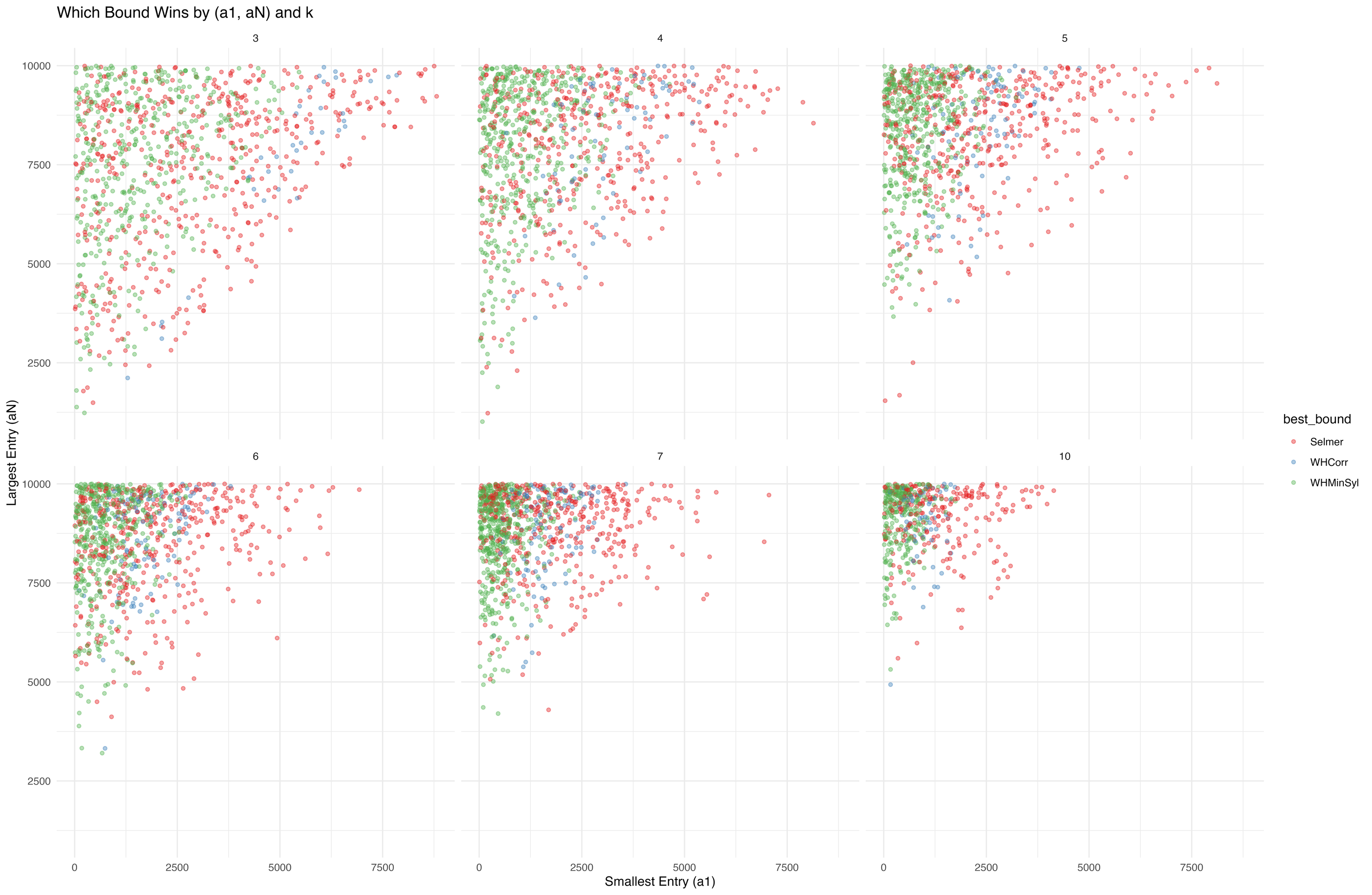}
\caption{Region-wise best-performing bound by \((a_1, a_n)\) with \(\|\boldsymbol{a}\|_\infty \le 10{,}000\).}
\label{fig:best_bound_regions_m10000}
\end{figure}

\newpage
\section{Proportion of Instances Where Each Bound is Tightest}
\label{appendix:best_bound_proportions}

The following figure summarises how often each of the considered bounds achieves the smallest error across randomly generated instances, broken down by dimension $n$. Note that each stacked bar shows the proportion of instances in which a given bound was tighter than the others.

\vspace{2.0mm}

\begin{figure}[ht!]
\centering
\includegraphics[width=1\linewidth]{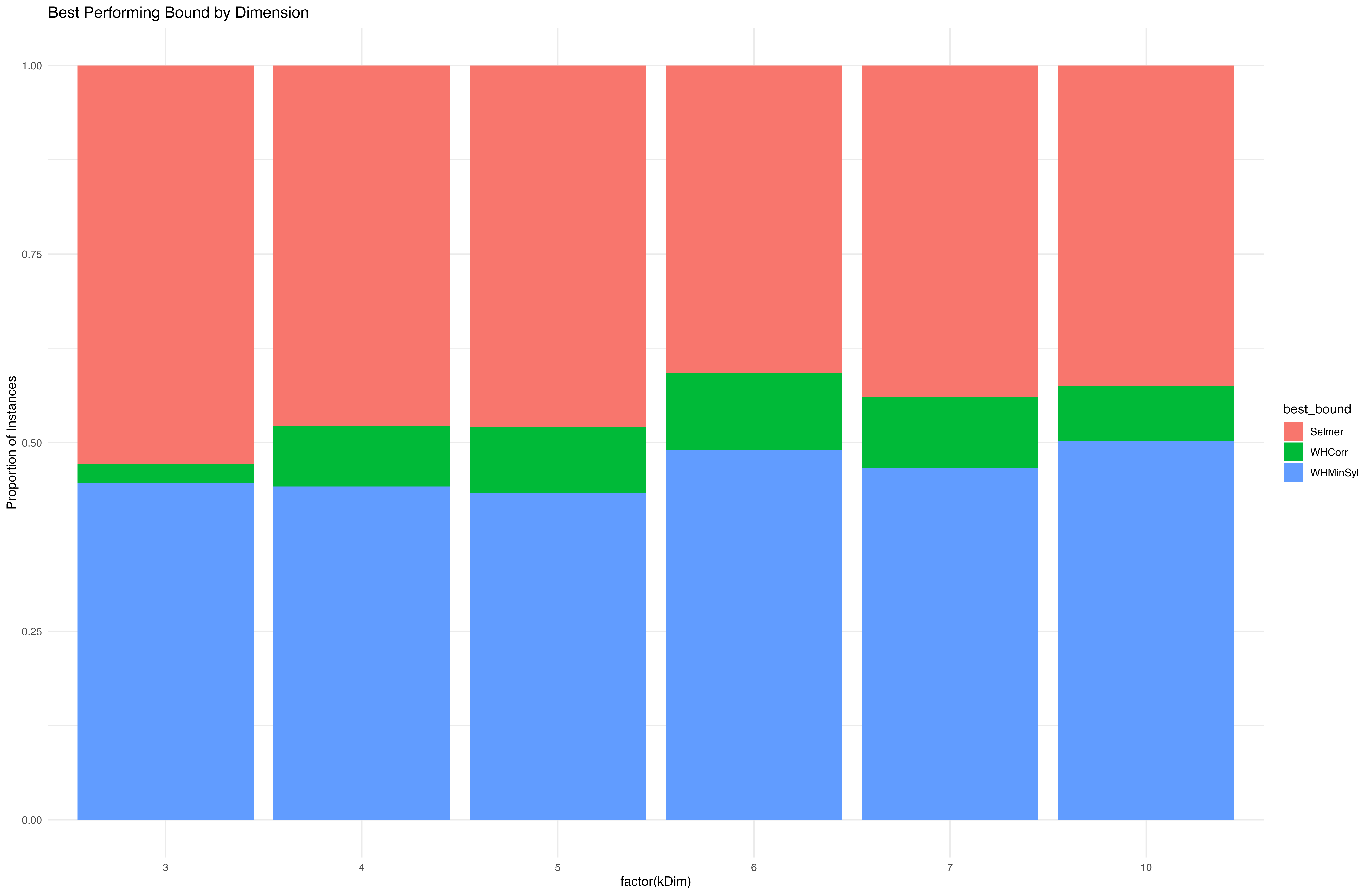}
\caption{Proportion of instances where each bound achieves the tightest bound on the Frobenius number $F(\boldsymbol{a})$, grouped by dimension $n$.}
\label{fig:stacked_bar_best_bound}
\end{figure}

\newpage
\section{Heat Map of the Errors Grouped by Dimension and Bound Type} 
\label{appendix:mean_error_heatmap_m10000}

The following heat map visualises the absolute error for each bound across bins of smallest and largest entry \((a_1, a_n)\), grouped by dimension. This complements Figure~\ref{fig:heatmap_all_bounds_structured} in the main text. Note that light colours (yellow) indicate larger overestimates, while darker colours (purple) indicate smaller errors, where the colour scale is log-transformed to improve visibility across the range of error values. 

\vspace{2.0mm}

\begin{figure}[ht!]
\centering
\includegraphics[width=\linewidth]{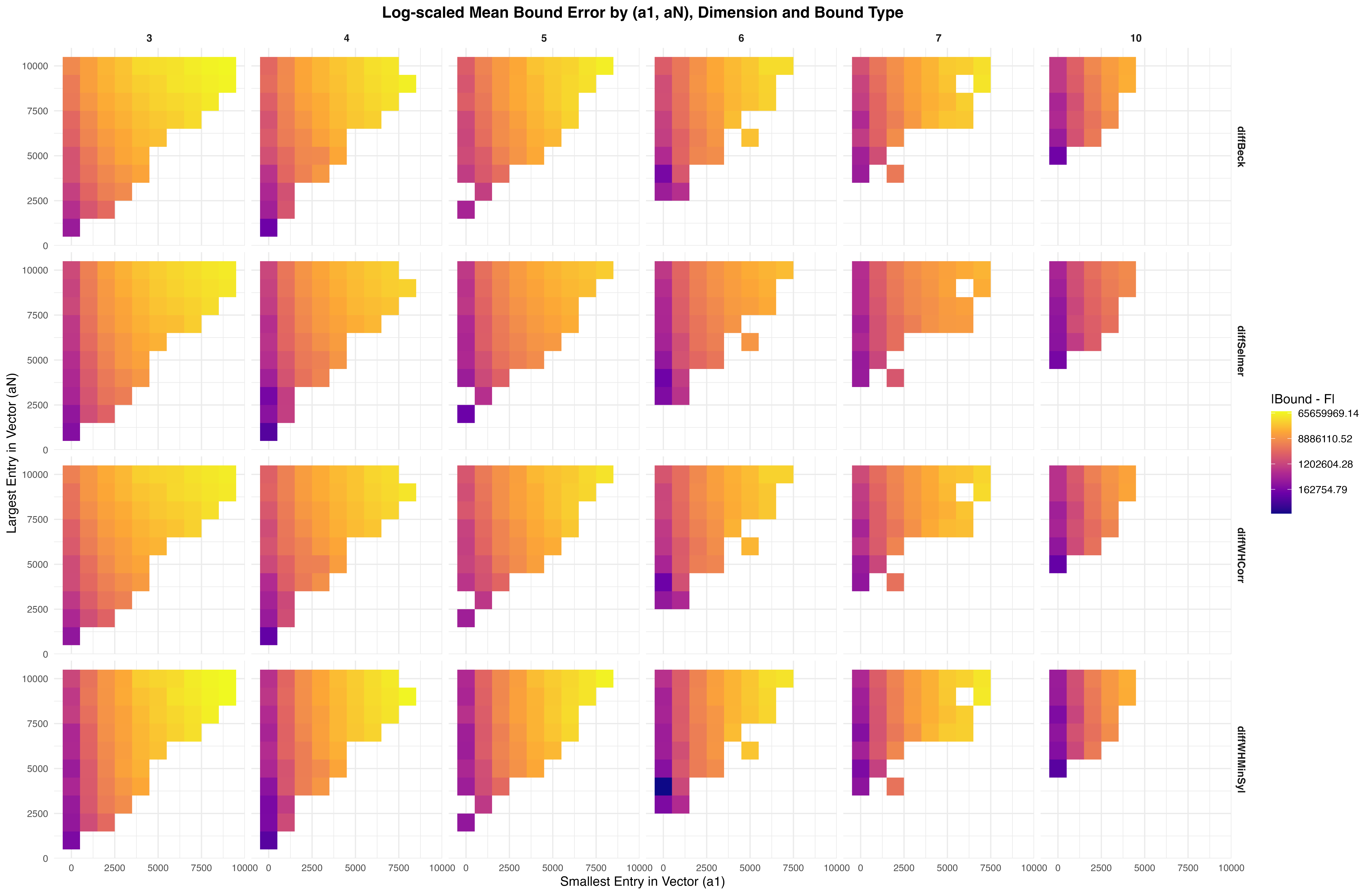}
\caption{Bound errors plotted by smallest and largest vector entries $(a_1, a_n)$ with
$\|\boldsymbol{a}\|_\infty \le 10{,}000$. Light colours (yellow) indicate larger overestimates, while darker colours (purple) indicate smaller errors. }
\label{fig:mean_bound_error_m10000_structured}
\end{figure}

\newpage
\section{Relative Error of the Williams and Haijima Bound \eqref{Williams_corrected_1} vs the Input Ratio $a_n/a_1$}\label{appendix:WH_Error_Input}
The following figure explores how structural imbalance in the input vector impacts upon the performance of the first bound of Williams and Haijima \eqref{Williams_corrected_1}. 
Here we consider the ratio \( a_n / a_1 \) as a proxy for vector conditioning, where a high ratio value indicates a large spread between the smallest and largest entries. Each point corresponds to a randomly generated instance. The horizontal axis corresponds to the ratio \( a_n / a_1 \) and the log-scaled vertical axis gives the relative error, namely the quantity 
$$
\frac{\texttt{WHCorr} - F(\boldsymbol{a})}{F(\boldsymbol{a})}.
$$

\vspace{2.0mm}

\begin{figure}[ht!]
\centering
\includegraphics[width=0.85\textwidth]{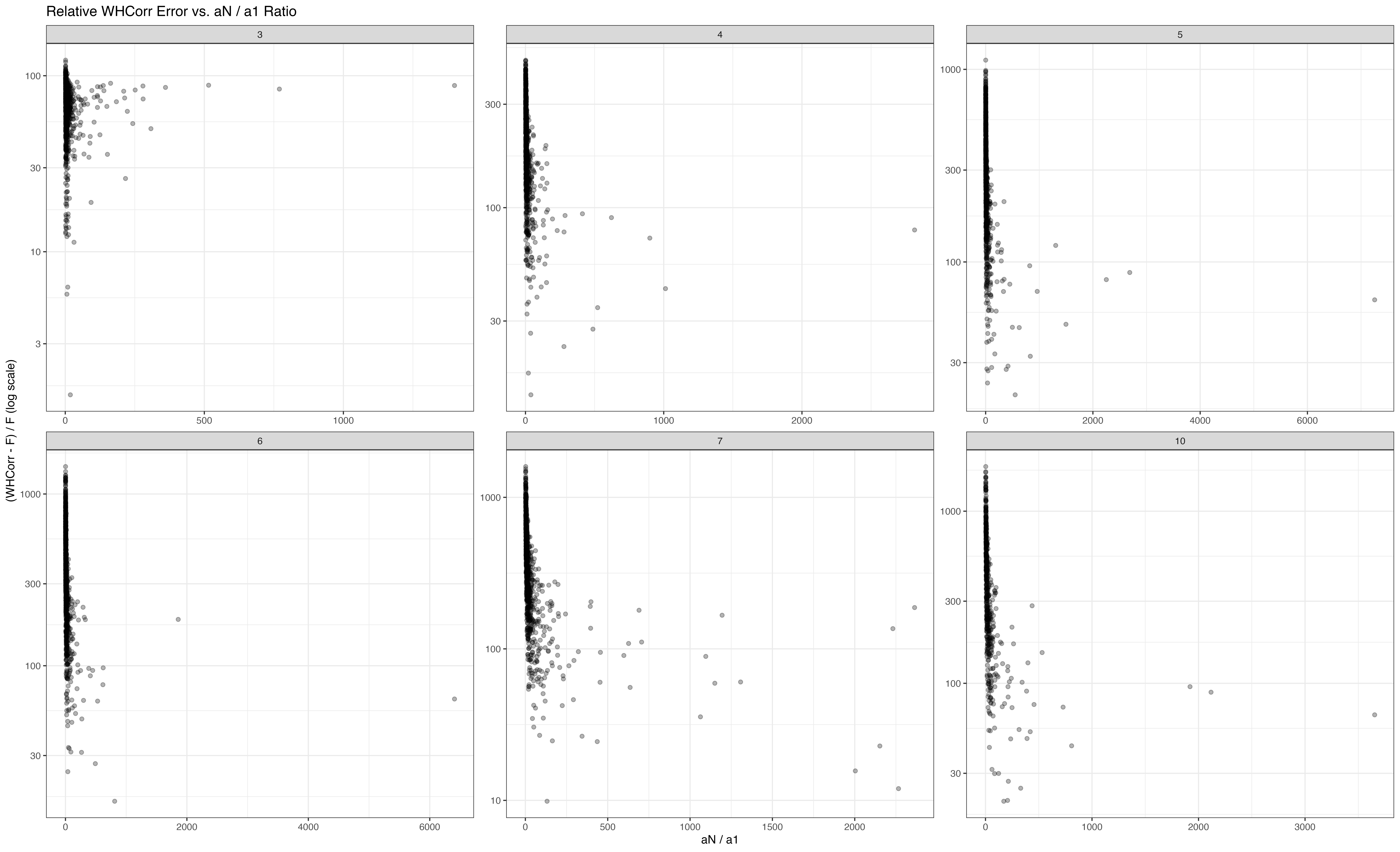}
\caption{The relative error of the first bound of Williams and Haijima \eqref{Williams_corrected_1} plotted against the input ratio $a_n / a_1$, grouped by dimension.}
\label{fig:rel_error_ratio_appendix}
\end{figure}


\end{document}